\documentclass[12pt, reqno]{amsart}

\usepackage[margin=1.2in]{geometry}
\usepackage{amsmath,amssymb,mathtools}
\usepackage{needspace}
\usepackage{tikz}
\usepackage{hyperref}
\usepackage{xcolor}

\DeclareMathOperator*{\argmin}{arg\,min}
\DeclareMathOperator{\tr}{tr}

\newtheorem{theorem}{Theorem}[section]
\newtheorem{proposition}[theorem]{Proposition}
\newtheorem{lemma}[theorem]{Lemma}
\newtheorem{corollary}[theorem]{Corollary}
\theoremstyle{definition}
\newtheorem{definition}[theorem]{Definition}
\theoremstyle{remark}
\newtheorem{remark}[theorem]{Remark}
\numberwithin{equation}{section}

\newcommand{\proofstep}[1]{\par\smallskip\Needspace{3\baselineskip}\noindent\emph{#1}\enspace}

\newcommand{\R}{\mathbb{R}}
\newcommand{\cB}{\mathcal{B}}

\newcommand{\cE}{\mathcal{E}}
\newcommand{\cN}{\mathcal{N}}
\newcommand{\cP}{\mathcal{P}}
\newcommand{\cS}{\mathcal{S}}
\newcommand{\cU}{\mathcal{U}}
\newcommand{\cW}{\mathcal{W}}
\newcommand{\al}{\alpha}
\newcommand{\bt}{\beta}
\newcommand{\ka}{\kappa}
\newcommand{\lam}{\lambda}
\newcommand{\Del}{\Delta}
\newcommand{\vphi}{\varphi}
\newcommand{\ep}{\varepsilon}
\newcommand{\grad}{\nabla}
\newcommand{\supp}{\operatorname{supp}}
\newcommand{\ibt}{\mathbf 1_{\{\bt=2\}}}

\definecolor{darkspringgreen}{rgb}{0.09, 0.45, 0.27} 
\definecolor{darkgray}{rgb}{0.66, 0.66, 0.66}
\definecolor{DeepBlue}{RGB}{0,0,160} 

\definecolor{darkblue}{rgb}{0.0,0.0,0.55}

\begin{document}

\title[Simplex ground states with small diffusion]{Uniqueness and shape of simplex ground states for aggregation--diffusion energies with mild repulsion}

\author{Tongseok Lim} \address{Tongseok Lim: Department of Quantitative Methods, Mitchell E. Daniels, Jr. School of Business\newline Purdue University, West Lafayette, IN 47907, USA} \email{lim336@purdue.edu}

\author{Geuntaek Seo}
\address{Geuntaek Seo: Department of Mathematics\newline POSTECH, Pohang Gyeongbuk 37673, Republic of Korea}
\email{gtseo@postech.ac.kr}



\begin{abstract}
We study global minimizers in $\R^d$ of the free energy $\mathcal E_\kappa=\mathcal W+\kappa\mathcal U$, where $\mathcal W$ is a power-law interaction energy with mild repulsion and $\mathcal U$ is the porous-medium internal energy with exponent $m>1$. In a range of exponents that we call the strict simplex regime, results of Lim--McCann and Davies--Lim--McCann show that $\mathcal W$ is minimized, uniquely up to rigid motions, by equal point masses at the vertices of a regular unit $d$-simplex. Diffusion excludes such atomic minimizers for every $\kappa>0$. We prove that uniqueness nevertheless persists: for every fixed $m>1$ and all sufficiently small $\kappa>0$, the free energy has a unique global minimizer up to rigid motions. After a suitable rigid motion, this minimizer is invariant under every symmetry of the simplex, and its support splits into $d+1$ clusters near the vertices, each of mass exactly $1/(d+1)$; in particular, it is not radially symmetric when $d\geq2$. The clusters have width of order $\kappa^{1/(d(m-1)+2)}$, which is also the exact order of the quadratic Wasserstein distance from the minimizer to the atomic simplex measure. Rescaled at this scale, the clusters converge in $L^1$ and $L^\infty$ to explicit Barenblatt-type profiles supported on ellipsoids determined by the Hessian of the simplex potential. For $d\geq2$ these ellipsoids are spheroids flattened along the directions of the vertices; when the repulsive exponent exceeds two, their aspect ratio equals $(d+1)^{-1/2}$, independently of all other parameters. At the endpoint $(4,2)$ of this regime in dimensions $d\geq2$, where the zero-diffusion minimizer is not unique, we show instead that vanishing diffusion selects the uniform spherical shell.
\end{abstract}

\maketitle

\section{Introduction}\label{section_intro}
For certain radial attractive--repulsive potentials, the interaction
energy is minimized, uniquely up to rigid motions, by equal point
masses at the vertices of a regular simplex \cite{LM21,DLM23}.  Adding
porous-medium diffusion rules out these atomic minimizers, since point
masses have infinite internal energy.  We show that rigidity
nevertheless survives: for sufficiently small diffusion, the free energy
has a unique global minimizer up to rigid motions.  Its density forms
separated clusters near the simplex vertices and, in dimensions two and
higher, is not radially symmetric.  Simplex symmetry and exact equality
of the cluster masses follow from uniqueness, without being imposed
in advance.  We also determine the concentration scale and the
limiting profiles of the clusters.

We consider the free energy
\begin{equation}\label{intro-free-energy}
 \cE_\ka(\rho)=\cW(\rho)+\ka\cU(\rho),\qquad\ka>0,
\end{equation}
on the Borel probability measures $\cP(\R^d)$, where
\[
 \cW(\rho):=\frac12\iint_{\R^d\times\R^d}
          W_{\al,\bt}(x-y)\,d\rho(x)\,d\rho(y),
\]
\[
 W_{\al,\bt}(x)=w_{\al,\bt}(|x|),  \quad w_{\al,\bt}(s)=s^\alpha/\alpha-s^\beta/\beta, \qquad\al>\bt\geq2,
\]
and
\[
 \cU(\rho):=
 \begin{cases}
 \displaystyle\frac1{m-1}\int_{\R^d}\rho(x)^m\,dx,
       &\rho\in\cP_{ac}(\R^d),\\[4pt]
 +\infty,&\text{otherwise},
 \end{cases}
 \qquad m>1.
\]
We call $\cW$ the interaction energy and $\cU$ the internal energy.
Here $\cP_{ac}(\R^d)$ denotes the absolutely continuous probability
measures.  We also write $W=W_{\al,\bt}$, $w=w_{\al,\bt}$, and
$U(s)=s^m/(m-1)$.  Since $w$ is minimized at $s=1$, the kernel is
repulsive at distances below one and attractive at distances above one.
The associated aggregation--diffusion equation is
\begin{equation}\label{pde}
 \partial_t\rho
 =\grad\cdot\big(\rho\grad(W*\rho)\big)+\ka\Delta\rho^m.
\end{equation}
Formally, this is the Wasserstein gradient flow of
\eqref{intro-free-energy}; see \cite{AGS08,San15}.  We study global
minimizers of the free energy, also called ground states; their role
for the long-time behavior of \eqref{pde} is studied in the companion
paper \cite{LS-dyn}.

Minimizers of interaction energies $\cW$ have been studied extensively.
Existence of global minimizers was established under sharp or almost
optimal conditions in \cite{CCP15,SST15}; see also \cite{CFT15} for
power-law potentials.  Balagu\'e, Carrillo, Laurent, and Raoul
\cite{BCLR13} related the dimension of the support of local minimizers
to the strength of the repulsion at the origin, and Carrillo, Delgadino,
and Mellet \cite{CDM16} showed, under suitable assumptions, that local
minimizers have bounded densities when the repulsion is at least as
singular as the Newtonian one.  Mild repulsion instead favors discrete
minimizers.  For a class of kernels with repulsive exponent strictly
larger than two, Carrillo, Figalli, and Patacchini \cite{CFP17} proved
that local minimizers in the $\infty$-Wasserstein topology have isolated
support points, and that global minimizers have finite support.  In one
dimension, Kang, Kim, Lim, and Seo \cite{KKLS21} studied uniqueness and
local minimality for power-law interactions, including how local
minimality depends on the Wasserstein topology.

Lim and McCann \cite{LM21} identified the global minimizer of $\cW$ for
each fixed $\bt\geq2$ and sufficiently large $\al$: it assigns equal
mass to the vertices of a regular unit $d$-simplex, up to rigid motions.
Their proof uses a sharp bound on the variance of measures with a fixed
support diameter.  For $\al>\bt>2$, they also showed that simplex
configurations with unequal masses can be local minimizers in the
$\infty$-Wasserstein topology.  Davies, Lim, and McCann \cite{DLM23}
then described the range of parameters for simplex minimizers through
a dimension-dependent threshold, above which the equally weighted
simplex is the unique global minimizer up to rigid motions.

Other shapes occur when $\bt=2$.  For $d\geq2$ and $2<\al<4$, Davies,
Lim, and McCann \cite{DLMshell22} proved that the unique global
minimizer up to translation is a spherical shell, that is, a uniform
measure on a sphere.  In these dimensions, at $(\al,\bt)=(4,2)$, there
are infinitely many centered minimizers: they are precisely the
probability measures supported on the sphere of radius
$\sqrt{d/(2d+2)}$ with zero first moment and second-moment matrix
$I_d/(2d+2)$ \cite{DLM23}.  For $d\geq2$, Frank and Matzke
\cite{FM25} extended the spherical-shell classification to a sharp range
of repulsive exponents $\bt\leq2$ for $2\leq\al\leq4$, with uniqueness
up to translations except at $(\al,\bt)=(4,2)$.  These classifications
concern the interaction energy without diffusion.  They show how
strongly the minimizing geometry depends on the exponents, but they do
not determine uniqueness once diffusion is added.  At the degenerate
pair $(4,2)$, a separate convexity argument gives uniqueness for
positive diffusion and selects the uniform shell as diffusion vanishes;
see Proposition~\ref{critical-selection}.

For energies with diffusion, Carrillo, Delgadino, and Patacchini
\cite{CDP19} studied existence and nonexistence of ground states and
derived their Euler--Lagrange conditions; see \cite{CCY19} for a survey
of aggregation--diffusion equations and their singular limits.  Burger
and Di Francesco \cite{BD08} constructed localized stationary profiles
for small quadratic diffusion in one dimension.  For smooth integrable
attractive kernels and quadratic diffusion, Burger, Di Francesco, and
Franek \cite{BDF13} characterized the existence of nontrivial steady
states and proved their uniqueness up to translation in one dimension,
using a strong form of the Krein--Rutman theorem.  Burger, Fetecau, and
Huang \cite{BFH14} studied compactly supported monotone equilibria in
one dimension for general power-law diffusion, including the different
free-boundary behavior on either side of the quadratic case.  For a
class of bounded smooth attractive potentials, Kaib \cite{Kai17} gave
existence criteria for compactly supported global minimizers for every
$m>1$ and proved uniqueness of stationary states up to translation in
the quadratic case whenever a minimizer exists.

Perturbative constructions also produce separated stationary clusters.
Fagioli and Jaafra \cite{FJ21} used a functional implicit-function
theorem for a one-dimensional two-species system with quadratic
diffusion.  They proved existence and uniqueness within a class of
multiple even bumps, with the individual bump masses and centers
prescribed, the centers satisfying a diffusion-free particle
equilibrium condition, and the effective well curvatures assumed
positive.  Our result concerns all global minimizers in arbitrary
dimension for every fixed $m>1$; the cluster masses and symmetry are
determined by minimization.

For purely attractive interactions, Carrillo, Hittmeir, Volzone, and Yao
\cite{CHVY19} proved radial symmetry of suitably regular stationary
states by continuous Steiner symmetrization.  Delgadino, Yan, and Yao
\cite{DYY22} then proved, under suitable assumptions, uniqueness up to
translation of bounded steady states with a fixed mass when $m\geq2$,
and constructed attractive potentials with nonunique steady states when
$1<m<2$.  The present kernels include short-range repulsion, and their
ground states have disconnected support and are nonradial when
$d\geq2$.  Thus the radial symmetry and interpolation arguments for
purely attractive interactions do not apply to these ground states.

The vanishing-diffusion limit has also been studied as a variational
approximation.  Fetecau, Kovacic, and Topaloglu \cite{FKT19} proved
convergence results for weak solutions and energy minimizers in domains
with boundaries, under their assumptions on the interaction and the
domain.  Craig and Topaloglu \cite{CT20} studied the different limit
$m\to\infty$, in which diffusion becomes a density constraint.

In this paper, $m > 1$ is fixed and $\ka\to0$, and our principal question is
uniqueness at each sufficiently small $\ka>0$, not merely the
identification of subsequential limits.  Compactness and energy
comparison identify the atomic limit, but they do not exclude distinct
diffuse minimizers converging to the same measure.  The strict
convexity of the internal energy does not settle this question, because
the interaction energy is nonconvex, and competing minimizers may a
priori have different cluster masses.  We therefore compare minimizers
without fixing these masses and prove uniqueness first; simplex
symmetry and exact equality of the masses then follow. 

We now fix the parameter range.  Let $d\geq1$, and let
$\al_\Delta(\bt)=\al_{\Delta^d}(\bt)$ be the transition threshold
of \cite[Theorem 1.5]{DLM23}.  We call
\begin{equation*}
 \cS:=\{(\al,\bt):\bt\geq2,\ \al>\al_\Delta(\bt)\}
\end{equation*}
the \emph{strict simplex regime}.  Here $\al_\Delta(\bt)\geq\bt$, and
when $\bt=2$ the strict condition is $\al>4$ for $d\geq2$ and $\al>3$
for $d=1$.  In this range, the potential generated by the simplex has
isolated minima at its vertices, with positive definite Hessians; see
Lemma~\ref{simplex-wells}.  A \emph{rigid motion} is a map
$x\mapsto Qx+a$, where $Q\in O(d)$ and $a\in\R^d$; reflections are
allowed.  We identify a simplex with its vertex set, and ``unit'' means
that each edge has length one.  We fix a centered regular unit simplex
and its equally weighted measure,
\[
 \Delta^d=\{x_0,\ldots,x_d\},
 \qquad
 \rho_\Delta=\frac{1}{d+1}\sum_{i=0}^d \delta_{x_i}.
\]
A \emph{symmetry} of $\Delta^d$ is a map $G\in O(d)$ that permutes
$\{x_0,\ldots,x_d\}$.

Our first theorem gives existence, uniqueness, symmetry, localization
near the simplex vertices, and exact equality of the cluster masses.

\begin{theorem}\label{main} Suppose $m>1$ and $(\al,\bt)\in\cS$, and put $r_\ka:=\ka^{1/(d(m-1)+2)}$. There exist $\ka_*>0$ and $L>0$ with $Lr_{\ka_*} \le 1/2$, depending only on $d$, $m$, $\al$, and $\bt$, such that for every $0<\ka<\ka_*$ the energy $\cE_\ka$ has a global minimizer, unique up to rigid motions. A representative is the normalized minimizer $\rho_\ka$ of Definition~\ref{normalized-definition}. Moreover, $G_\#\rho_\ka=\rho_\ka$ for every symmetry $G$ of $\Delta^d$, so $\rho_\ka$ is centered, and \[ \supp\rho_\ka\subset\bigcup_{i=0}^dB_{Lr_\ka}(x_i),\qquad \rho_\ka\big(B_{Lr_\ka}(x_i)\big)=\frac1{d+1},\quad i=0,\ldots,d. \] 
\end{theorem}

We impose no symmetry a priori and do not prescribe the individual cluster masses. For $d\geq2$ the minimizer $\rho_\ka$ is not radially symmetric, so the rearrangement and symmetrization arguments available for purely attractive kernels \cite{CHVY19,DYY22} are not applicable. To the best of our knowledge, Theorem~\ref{main} is the first result establishing uniqueness up to rigid motions of global minimizers of an aggregation--diffusion free energy on $\R^d$ with a radial interaction kernel and no external potential, in a setting where these minimizers are not radially symmetric. The simplex symmetry is a consequence of uniqueness; since the symmetries of $\Delta^d$ act transitively on its vertices, all cluster masses are equal. This equality holds for every sufficiently small $\ka>0$; convergence to the atomic simplex alone would yield equal masses only in the limit.

We also determine the size and shape of the clusters.  The scale
$r_\ka$ in Theorem~\ref{main} balances the interaction cost $r^2$ with
the diffusion cost $\ka r^{-d(m-1)}$.  Put $u:=1/(d+1)$,
$\hat x_i:=x_i/|x_i|$, and
\begin{equation}\label{well-constants}
 a_\perp:=\frac{\al-\bt}{2}-\ibt,
 \qquad
 a_\parallel:=\frac{(\al-\bt)(d+1)}{2}-\ibt.
\end{equation}

\begin{theorem}\label{shape}
Suppose $m>1$ and $(\al,\bt)\in\cS$, and let $\ka_*,L$ be as in Theorem~\ref{main}. For $\ka \in (0, \ka_*)$, let $\rho_\ka$ be the normalized minimizer of
Theorem~\ref{main}, with clusters
$\rho_{\ka,i}:=\rho_\ka\mathbf 1_{B_{1/2}(x_i)}$. Then $\|\rho_{\ka,i}\|_{L^\infty}\leq Lr_\ka^{-d}$, and each set $\{\rho_{\ka,i}>0\}$ is convex with $C^2$ boundary, with $\rho_\ka^{m-1}$ concave on it. Moreover, $r_\ka^d\rho_{\ka,i}(x_i+r_\ka z)\longrightarrow f_i(z)$ as $\ka \to 0$ in $L^1(\R^d)\cap L^\infty(\R^d)$, 
where
\begin{equation}\label{intro-profile}
 f_i(z)=\left[\frac{m-1}{m}
       \big(\lam_*-z\cdot A_i z\big)\right]_+^{1/(m-1)},
\end{equation}
\begin{equation}\label{intro-matrix}
 A_i:=\frac12D^2(W*\rho_\Delta)(x_i)
 =\frac u2\Big(a_\perp\big(I_d-\hat x_i\otimes\hat x_i\big)
     +a_\parallel\,\hat x_i\otimes\hat x_i\Big),
\end{equation}
and $\lam_*>0$ is the unique constant satisfying $\int f_i=u$ for every $i$, given explicitly in \eqref{lambda-explicit}. 
The rescaled supports
converge in Hausdorff distance to the ellipsoids
$E_i:=\{z\in\R^d:z\cdot A_iz\leq\lam_*\}$. Additionally, the quadratic Wasserstein distance satisfies
\[
 \lim_{\ka\to0}
 \frac{d_2(\rho_\ka,\rho_\Delta)^2}{r_\ka^2}
 =
 \frac{2(d+1)\lam_*}{d+2+2/(m-1)}
 \left(\frac{1}{a_\parallel}+\frac{d-1}{a_\perp}\right),
\]
where the term involving $a_\perp$ is omitted when $d=1$.
\end{theorem}

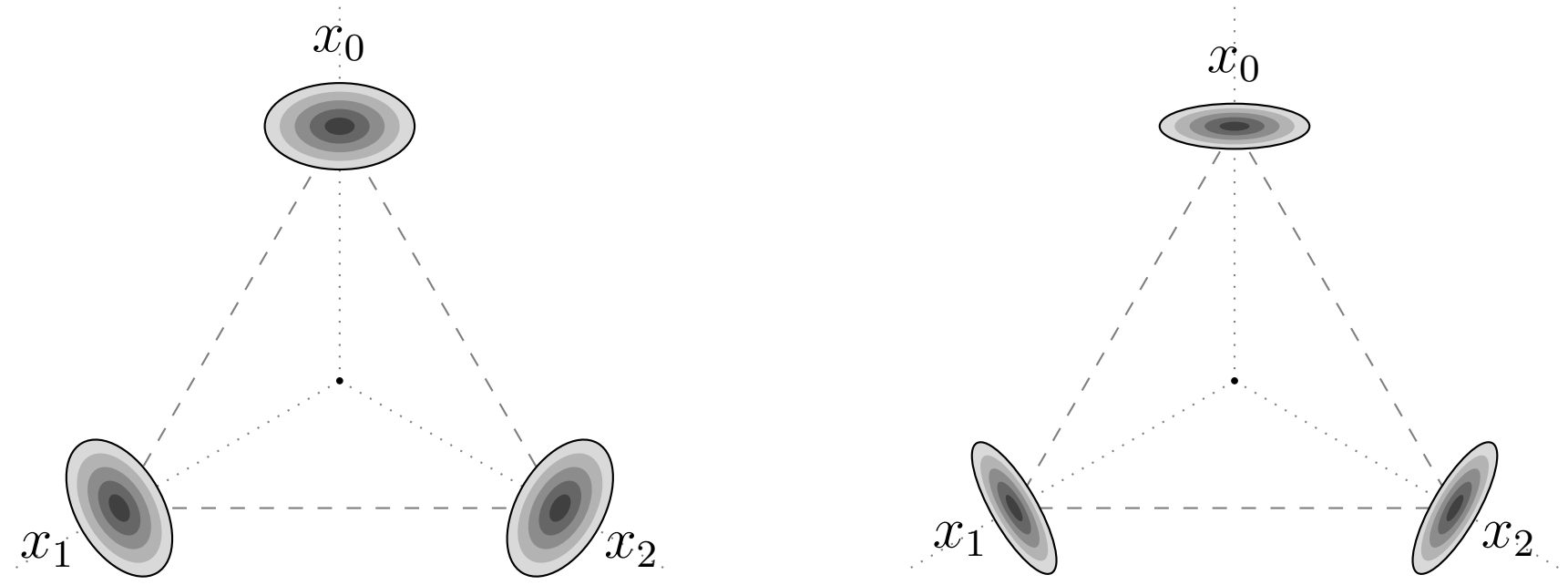
\begin{figure}[t]
\centering
\begin{tikzpicture}[scale=3.2]
\draw[dashed,gray] (0.0000,0.5774) -- (-0.5000,-0.2887) -- (0.5000,-0.2887) -- cycle;
\draw[dotted,gray] (0,0) -- (0.0000,0.8600);
\draw[dotted,gray] (0,0) -- (-0.7448,-0.4300);
\draw[dotted,gray] (0,0) -- (0.7448,-0.4300);
\fill (0,0) circle (0.008);
\begin{scope}[shift={(0.0000,0.5774)},rotate=90]
  \fill[black!15] (0,0) ellipse (0.0981 and 0.1700);
  \fill[black!30] (0,0) ellipse (0.0785 and 0.1360);
  \fill[black!45] (0,0) ellipse (0.0589 and 0.1020);
  \fill[black!60] (0,0) ellipse (0.0393 and 0.0680);
  \fill[black!75] (0,0) ellipse (0.0196 and 0.0340);
  \draw[thin] (0,0) ellipse (0.0981 and 0.1700);
\end{scope}
\node at (0.0000,0.7655) {$x_0$};
\begin{scope}[shift={(-0.5000,-0.2887)},rotate=210]
  \fill[black!15] (0,0) ellipse (0.0981 and 0.1700);
  \fill[black!30] (0,0) ellipse (0.0785 and 0.1360);
  \fill[black!45] (0,0) ellipse (0.0589 and 0.1020);
  \fill[black!60] (0,0) ellipse (0.0393 and 0.0680);
  \fill[black!75] (0,0) ellipse (0.0196 and 0.0340);
  \draw[thin] (0,0) ellipse (0.0981 and 0.1700);
\end{scope}
\node at (-0.6629,-0.3827) {$x_1$};
\begin{scope}[shift={(0.5000,-0.2887)},rotate=330]
  \fill[black!15] (0,0) ellipse (0.0981 and 0.1700);
  \fill[black!30] (0,0) ellipse (0.0785 and 0.1360);
  \fill[black!45] (0,0) ellipse (0.0589 and 0.1020);
  \fill[black!60] (0,0) ellipse (0.0393 and 0.0680);
  \fill[black!75] (0,0) ellipse (0.0196 and 0.0340);
  \draw[thin] (0,0) ellipse (0.0981 and 0.1700);
\end{scope}
\node at (0.6629,-0.3827) {$x_2$};
\node[below] at (0,-0.48) {\footnotesize (a) $\beta>2$: ratio $1/\sqrt{3}\approx0.58$};
\end{tikzpicture}
\hspace{0.6cm}
\begin{tikzpicture}[scale=3.2]
\draw[dashed,gray] (0.0000,0.5774) -- (-0.5000,-0.2887) -- (0.5000,-0.2887) -- cycle;
\draw[dotted,gray] (0,0) -- (0.0000,0.8600);
\draw[dotted,gray] (0,0) -- (-0.7448,-0.4300);
\draw[dotted,gray] (0,0) -- (0.7448,-0.4300);
\fill (0,0) circle (0.008);
\begin{scope}[shift={(0.0000,0.5774)},rotate=90]
  \fill[black!15] (0,0) ellipse (0.0513 and 0.1700);
  \fill[black!30] (0,0) ellipse (0.0410 and 0.1360);
  \fill[black!45] (0,0) ellipse (0.0308 and 0.1020);
  \fill[black!60] (0,0) ellipse (0.0205 and 0.0680);
  \fill[black!75] (0,0) ellipse (0.0103 and 0.0340);
  \draw[thin] (0,0) ellipse (0.0513 and 0.1700);
\end{scope}
\node at (0.0000,0.7186) {$x_0$};
\begin{scope}[shift={(-0.5000,-0.2887)},rotate=210]
  \fill[black!15] (0,0) ellipse (0.0513 and 0.1700);
  \fill[black!30] (0,0) ellipse (0.0410 and 0.1360);
  \fill[black!45] (0,0) ellipse (0.0308 and 0.1020);
  \fill[black!60] (0,0) ellipse (0.0205 and 0.0680);
  \fill[black!75] (0,0) ellipse (0.0103 and 0.0340);
  \draw[thin] (0,0) ellipse (0.0513 and 0.1700);
\end{scope}
\node at (-0.6223,-0.3593) {$x_1$};
\begin{scope}[shift={(0.5000,-0.2887)},rotate=330]
  \fill[black!15] (0,0) ellipse (0.0513 and 0.1700);
  \fill[black!30] (0,0) ellipse (0.0410 and 0.1360);
  \fill[black!45] (0,0) ellipse (0.0308 and 0.1020);
  \fill[black!60] (0,0) ellipse (0.0205 and 0.0680);
  \fill[black!75] (0,0) ellipse (0.0103 and 0.0340);
  \draw[thin] (0,0) ellipse (0.0513 and 0.1700);
\end{scope}
\node at (0.6223,-0.3593) {$x_2$};
\node[below] at (0,-0.48) {\footnotesize (b) $\beta=2,\ \alpha=4.5$: ratio $1/\sqrt{11}\approx0.30$};
\end{tikzpicture}
\caption{Limiting cluster shapes of Theorem~\ref{shape} for $d=2$.  The
dashed triangle is the unit simplex $\Delta^2$ and the dotted rays are
the axes $\hat x_i$.  Each cluster is drawn as the level sets of its
limiting profile $f_i$, homothetic ellipses centered at the vertex.  The
cluster widths, of order $r_\ka$, are magnified by a common factor;
the shapes and orientations are exact.  Each ellipse is flattened along
$\hat x_i$ with semi-axis ratio $\sqrt{a_\perp/a_\parallel}$: (a) for every $(\al,\bt) \in \cS$ with $\bt>2$ and $m>1$, this
ratio is $1/\sqrt3$; (b) for $\bt=2$ it depends on $\al$, here
$\al=4.5$, and it tends to $0$ as $\al\downarrow4$.}
\label{fig-clusters}
\end{figure}

In the strict simplex regime, $a_\parallel>0$, and $a_\perp>0$ when
$d\geq2$; see Lemma~\ref{simplex-wells}.  Thus each limiting profile is
a Barenblatt-type profile on the ellipsoid $E_i$.  For $d\geq2$, $E_i$
is a spheroid with symmetry axis $\hat x_i$, flattened along this axis:
the ratio of its semi-axis in the direction $\hat x_i$ to its other
semi-axes is
\[
 \sqrt{\frac{a_\perp}{a_\parallel}}
 =\begin{cases}
 (d+1)^{-1/2},&\bt>2,\\[6pt]
 \displaystyle\sqrt{\frac{\al-4}{(\al-2)(d+1)-2}},&\bt=2.
 \end{cases}
\]
In particular, these limiting ellipsoids are not balls although the
pair potential is radial.  The cutoff power in \eqref{intro-profile} is
the classical Barenblatt form of porous-medium diffusion \cite{V07},
while $A_i$ is half the Hessian of the simplex potential
$W*\rho_\Delta$ at $x_i$.  Its anisotropic part comes from the other
clusters, and this anisotropic confinement makes the clusters
non-spherical.  When $\bt>2$, the aspect ratio does not depend on $\al$,
$\bt$, or $m$.  When $\bt=2$, it tends to zero as $\al\downarrow4$,
where the simplex regime meets the spherical-shell regime of
\cite{DLMshell22}.  By Theorem~\ref{shape}, $r_\ka$ is also the exact
order of the quadratic Wasserstein distance to $\rho_\Delta$.
Figure~\ref{fig-clusters} illustrates the two cases for $d=2$.

\begin{remark}[Dynamical stability]\label{dynamics-remark}
Transport and entropy methods for quantitative convergence to
equilibrium were developed, in particular, in \cite{CMV03,CMV06}.  The
relation between these ground states and the long-time behavior of
\eqref{pde} is studied in the companion paper \cite{LS-dyn}, where the
algebraic identity of Lemma~\ref{dipole-positivity} enters a local
displacement-convexity argument that gives exponential convergence to a
fixed rigid image of $\rho_\ka$ for nearby finite-energy data in the
$\infty$-Wasserstein topology, with rates uniform in small $\ka$.  The
static results proved here do not depend on that analysis.  For
comparison, Di Francesco and Jaafra \cite{DFJ19} established local
asymptotic stability for a one-dimensional equation with quadratic
diffusion and a purely attractive integrable kernel.
\end{remark}

\subsection*{Why the small-diffusion limit is singular}
The atomic minimizer $\rho_\Delta$ lies outside the finite-energy
domain of $\cU$ for every $\ka>0$, so a smooth perturbation argument
cannot be applied directly to $\cE_\ka$ at $\rho_\Delta$.  This does not
preclude implicit-function constructions in rescaled variables, as in
\cite{FJ21}; the difficulty here is to compare arbitrary minimizers
while their cluster masses, supports, and symmetries are still unknown.
After rescaling, the free boundaries move with the minimizer, and the
diffusion weight $U''(f)=mf^{m-2}$ vanishes there when $m>2$ and is
singular when $1<m<2$.  Rigid motions introduce further degeneracies.
We treat these issues by comparing two minimizers directly through
their first variations, after an explicit alignment that removes
translations and rotations; see Section~\ref{subsection-uniqueness}.

The proof of Theorem~\ref{main} starts with a competitor obtained by
smoothing the simplex atoms, which gives an $O(r_\ka^2)$ bound on the
energy excess above $\cW(\rho_\Delta)$.  The global Euler--Lagrange
equation has one common Lagrange multiplier even though the support is
disconnected; it confines the clusters to the non-degenerate wells of
the simplex potential and determines their width.  We then remove the
rigid-motion freedom by an alignment condition that involves only the
masses and first moments of the clusters.

Two estimates close the comparison of aligned minimizers.  First, after
the differences of cluster masses have been separated off, the limiting
interaction depends only on the first moments, and the simplex geometry
makes this quadratic form nonnegative under alignment
(Lemma~\ref{dipole-positivity}).  Second, the concavity of the rescaled
pressure gives a uniform bound on the inverse diffusion weights up to
the free boundaries (Lemma~\ref{inverse-weight-lemma}).  A weighted
Cauchy--Schwarz estimate then bounds the $L^1$ difference of the
rescaled densities by their diffusion contribution, so the interaction
errors can be absorbed for small $\ka$, while a positive zeroth-order
term controls the differences of the cluster masses.  Thus the argument
proves uniqueness without assuming equal masses, for every fixed $m>1$,
including the free-boundary degeneracy when $m>2$.

For Theorem~\ref{shape}, we rescale the normalized minimizers, which
are $\Delta^d$-symmetric by Theorem~\ref{main}, at the scale $r_\ka$.  Translating
one cluster at a time gives force-balance equations, which together
with the symmetry force the limiting cluster means to vanish.  The
Euler--Lagrange equation then identifies the profiles.

Section~\ref{section-preliminaries} gives the notation and variational
tools.  Section~\ref{section-uniqueness} proves Theorem~\ref{main}, and
Section~\ref{section-shape} proves Theorem~\ref{shape}.
Section~\ref{section-remarks} contains concluding remarks, including the
selection result at the degenerate pair $(4,2)$ and open questions.

\section{Preliminaries}\label{section-preliminaries}
\subsection{Notation and elementary bounds}
We write $\cP(\R^d)$ for Borel probability measures and $\cP_{ac}(\R^d)$
for those that are absolutely continuous with respect to Lebesgue measure.
For $p\geq1$, $\cP_p(\R^d)$ denotes the measures with finite $p$-th moment,
and $d_p$ is the $p$-Wasserstein distance.  In contrast, the subscript in
$\cP_0(\R^d)$ denotes centering:
\[
 \cP_0(\R^d):=\{\rho\in\cP_1(\R^d):\textstyle\int x\,d\rho(x)=0\}.
\]
Narrow convergence means convergence against bounded continuous
functions, and $T_\#\rho$ denotes the push-forward.  For a finite
non-negative measure, $|\rho|$ denotes its total mass.  We take $B_r(x)$
to be the open Euclidean ball, put $B_r:=B_r(0)$, write
$[s]_+=\max\{s,0\}$, and use $\supp\rho$ for the closed support.  The
same symbol denotes an absolutely continuous measure and its density.
Unqualified integrals are over $\R^d$.  For vectors,
$(v\otimes w)_{kl}=v_kw_l$.  For matrices,
$\langle M,N\rangle:=\tr(M^TN)$ and $|M|:=\langle M,M\rangle^{1/2}$ denote the Frobenius inner
product and norm, and $\|M\|:=\sup_{|x|\leq1}|Mx|$ denotes the operator norm.

Let $\Delta^d=\{x_0,\ldots,x_d\}$ be the reference simplex fixed in the
introduction, and set
\[
 u:=\frac1{d+1},\qquad \rho_\Delta=u\sum_{i=0}^d\delta_{x_i},\qquad
 \cP_\Delta:=\cP_0(\R^d)\cap\argmin\cW,\qquad e_{ij}:=x_i-x_j.
\]
For $(\al,\bt)\in\cS$, \cite[Theorem 1.5]{DLM23} shows that
$\cP_\Delta=\{Q_\#\rho_\Delta:Q\in O(d)\}$.  Up to an orthogonal map, $x_j=(\mathbf e_j-\bar{\mathbf e})/\sqrt2$,
where $\mathbf e_0,\ldots,\mathbf e_d$ is the standard basis of
$\R^{d+1}$, $\bar{\mathbf e}:=u\sum_k\mathbf e_k$, and the hyperplane
$\{y\in\R^{d+1}:\sum_ky_k=0\}$ is identified with $\R^d$.  For $v$ in
this hyperplane, $\sum_j(x_j\cdot v)^2=|v|^2/2$.  This gives
\begin{equation}\label{simplex-geometry}
 \begin{gathered}
 \sum_{j=0}^dx_j=0,\qquad
 |x_i|^2=\frac{d}{2(d+1)},\qquad
 x_i\cdot x_j=-\frac1{2(d+1)}\quad(i\neq j),\\
 \sum_{j=0}^dx_j\otimes x_j=\frac12I_d.
 \end{gathered}
\end{equation}
In particular $|e_{ij}|=1$ for $i\neq j$, and the balls $B_{1/2}(x_i)$
are pairwise disjoint.  For a measure $\rho$ supported in
$\bigcup_iB_{1/2}(x_i)$, the restrictions
$\rho_i:=\rho\mathbf 1_{B_{1/2}(x_i)}$ are called its \emph{clusters}.
For a set of measures $\mathcal A$, we write
$d_\al(\rho,\mathcal A)=\inf_{\nu\in\mathcal A}d_\al(\rho,\nu)$.
Constants denoted by $c,C$ may change from line to line; unless stated
otherwise they depend only on $d$, $m$, $\al$, and $\bt$, not on a
sufficiently small $\ka$ or on the choice of global minimizer.

The minimum value of $W$ is $w(1)=1/\al-1/\bt<0$, and we put $c_W:=-w(1)>0$.  There are constants
$c_0,C_0>0$ such that
\begin{equation}\label{polynomial-coercivity}
 c_0|z|^\al-C_0\leq W(z),\qquad
 |W(z)|\leq C_0(1+|z|^\al).
\end{equation}
In particular, finite interaction energy implies a finite $\al$-th moment:
Fubini first gives a point $y_0$ with $\int|x-y_0|^\al\,d\rho(x)<\infty$.
If $b=\int x\,d\rho(x)$, Jensen's inequality then gives
\begin{equation}\label{centered-moment-bound}
 \int|x-b|^\al\,d\rho(x)
 \leq\iint|x-y|^\al\,d\rho(x)d\rho(y)
 \leq\frac{2\cW(\rho)+C_0}{c_0}.
\end{equation}
Also $W\in C^2(\R^d)$, and its derivatives of order at most two have at
most $\al$-power growth.  Thus $W*\rho$ is $C^2$ whenever
$\rho\in\cP_\al(\R^d)$, by differentiation under the integral.  These observations will justify both centering and the continuous
Euler--Lagrange representative before a support bound is known.

\subsection{The global Euler--Lagrange equation}
A global minimizer has one common Lagrange multiplier, even when its
support is disconnected.  This will allow us to compare cluster masses.
We give the first-variation argument; see also \cite{CDP19,San15}.

\begin{proposition}\label{CDP_prop}
Let $\rho$ be a global minimizer of $\cE_\ka$, and set
\[
 \Phi:=W*\rho,\qquad \ell:=\ka m\cU(\rho)+2\cW(\rho).
\]
Then
\begin{equation}\label{EL_global}
 \begin{cases}
 \ka U'(\rho)+\Phi=\ell,&\rho\text{-almost everywhere},\\
 \ka U'(\rho)+\Phi\geq\ell,&\text{Lebesgue-almost everywhere}.
 \end{cases}
\end{equation}
\end{proposition}
\begin{proof}
The identities $\int U'(\rho)\rho=m\cU(\rho)$ and
$\int\Phi\,d\rho=2\cW(\rho)$ give
$\ell=\int(\ka U'(\rho)+\Phi)\,d\rho$.
For any bounded compactly supported probability density $\nu$,
minimality along $(1-t)\rho+t\nu$ gives
\begin{equation}\label{first_variation}
 0\leq\left.\frac{d}{dt}\right|_{t=0+}
       \cE_\ka((1-t)\rho+t\nu)
   =\int(\ka U'(\rho)+\Phi)(\nu-\rho).
\end{equation}
The derivative is justified by $\rho\in L^m$, the growth of $U'$,
and \eqref{polynomial-coercivity}.  Testing with normalized indicators
of bounded measurable sets shows that $\ka U'(\rho)+\Phi\geq\ell$
almost everywhere.  Its integral against $\rho$ equals $\ell$, so
equality holds $\rho$-almost everywhere.
\end{proof}

\begin{corollary}[The continuous density]\label{version}
With $\ell$ as above, a global minimizer $\rho$ has the continuous
representative
\begin{equation}\label{rho-version}
 \rho(x)=\left[\frac{m-1}{m\ka}
           \big(\ell-W*\rho(x)\big)\right]_+^{1/(m-1)}.
\end{equation}
\end{corollary}
\begin{proof}
On $\{\rho>0\}$, the equality in \eqref{EL_global} gives this formula
almost everywhere.  On $\{\rho=0\}$, the inequality gives
$W*\rho\geq\ell$, so the right-hand side is zero.  It is continuous
because $W*\rho$ is $C^2$.
\end{proof}
We use this representative from now on.  For a global minimizer
$\rho_\ka$ we write $\ell_\ka:=\ka m\cU(\rho_\ka)+2\cW(\rho_\ka)$ and
$\Phi_\ka:=W*\rho_\ka$, so that
\[
 \rho_\ka(x)>0\quad\Longleftrightarrow\quad\Phi_\ka(x)<\ell_\ka,
 \qquad
 \supp\rho_\ka\subset\{\Phi_\ka\leq\ell_\ka\}.
\]

\subsection{Compactness up to translations}
For completeness we give the compactness argument for our power-law
kernel.  It is a special case of the coercive-interaction compactness
principle in \cite[Section 2]{CDP19}.  

\begin{lemma}\label{CDP_lemma}
If $(\rho_n)\subset\cP(\R^d)$ and $\liminf_n\cW(\rho_n)<\infty$,
then a subsequence converges narrowly to a probability measure after
translations.
\end{lemma}
\begin{proof}
Choose a subsequence with uniformly bounded interaction energy, and
translate each measure by its barycenter.  The moments in
\eqref{centered-moment-bound} are uniformly bounded.  Consequently the
translated measures assign at most $CR^{-\al}$ mass to $B_R^c$, uniformly
along the subsequence.  They are tight, and Prokhorov's theorem gives the
claimed narrow subsequential convergence to a probability measure.
\end{proof}

\section{Proof of Theorem \ref*{main}: Uniqueness, symmetry, and equal masses}\label{section-uniqueness}
\subsection{Existence and zero-diffusion compactness}\label{subsection-compactness}
This subsection collects three facts about global minimizers as $\ka\to0$ used throughout the proof of Theorem~\ref{main}. Lemma~\ref{existence} gives existence and confines all centered minimizers to a fixed ball. Lemma~\ref{compactsupport} compares the minimizers with a smoothed copy of the atomic simplex measure. This gives an energy excess of order $r_\ka^2$, a bound of the same order for the Lagrange multipliers, and convergence of centered minimizers to the set of atomic minimizers $\cP_\Delta$. Lemma \ref{potential-convergence} turns this convergence of measures into $C^2$ convergence of the interaction potentials, which will be used to localize the clusters in the next subsection.

\begin{lemma}[Existence and uniform support]\label{existence}
For every $\ka>0$, the energy $\cE_\ka$ has a global minimizer.
There is $R>0$ such that every centered global minimizer with
$0<\ka\leq1$ is supported in $B_R$.
\end{lemma}
\begin{proof}
Let $(\rho_n)$ be a minimizing sequence.  Since
$\cE_\ka\geq\cW\geq w(1)/2$ and $\cU\geq0$, both $\cW(\rho_n)$ and
$\cU(\rho_n)$ are bounded uniformly in $n$.  By
Lemma~\ref{CDP_lemma}, after translations and passage to a subsequence,
$\rho_n$ converges narrowly to some $\rho\in\cP(\R^d)$.  Both $\cU$ and
$\cW$ are translation invariant, so the translated sequence is still
minimizing.  Both are also narrowly lower semicontinuous.  For $\cU$,
this holds because $\cU$ is the integral functional of the convex
superlinear function $U$ with $U(0)=0$; alternatively, since
$\cU(\rho_n)$ is bounded, a further subsequence converges weakly in
$L^m$, necessarily to $\rho$, and the $L^m$ norm is weakly lower
semicontinuous.  For $\cW$, the measures $\rho_n\otimes\rho_n$ converge
narrowly to $\rho\otimes\rho$; we apply this to the bounded continuous
truncations $\min\{W-w(1),k\}$ and let $k\to\infty$ by monotone
convergence, using that $W$ is continuous and bounded below.
Hence $\rho$ is a minimizer.  Since $\cE_\ka(\rho)<\infty$, $\rho \in L^m$, and $\rho$ has finite $\al$-th moment by
\eqref{centered-moment-bound}.

For the support bound, fix a smooth compactly supported probability
density $\nu$, and let $\rho_\ka$ denote a centered minimizer for $\cE_\ka$.  Minimality gives, for $0<\ka\leq1$,
\[
 \cW(\rho_\ka)\leq\cE_\ka(\rho_\ka)\leq\cE_1(\nu),
 \qquad
 \ka\cU(\rho_\ka)\leq\cE_1(\nu)-\frac{w(1)}2.
\]
Thus the Euler--Lagrange multiplier
$\ell_\ka=\ka m\cU(\rho_\ka)+2\cW(\rho_\ka)$ 
is bounded above independently of the minimizer and of $\ka\leq1$.
Jensen's inequality and
\eqref{polynomial-coercivity} give
\[
 \Phi_\ka(x) = W*\rho_\kappa (x)
 \geq c_0\int|x-y|^\al\,d\rho_\ka(y)-C_0
 \geq c_0|x|^\al-C_0.
\]
Choose $R$ so that this lower bound exceeds the multiplier bound for
$|x|\geq R$.  Formula \eqref{rho-version} then gives
$\supp\rho_\ka\subset B_R$.
\end{proof}

Set
\begin{equation*}
 r_\ka:=\ka^{1/(d(m-1)+2)}.
\end{equation*}
This scale balances the diffusion cost and the quadratic interaction
cost:
\begin{equation}\label{cluster-scale-balance}
 \ka r_\ka^{-d(m-1)}=r_\ka^2.
\end{equation}

\begin{lemma}[A simplex competitor and compactness]\label{compactsupport}
For $(\al,\bt)\in\cS$ and every sufficiently small $\ka>0$, every
global minimizer satisfies
\begin{equation}\label{energy-rate}
 0\leq\cE_\ka(\rho_\ka)-\cW(\rho_\Delta)
 =\cW(\rho_\ka)-\cW(\rho_\Delta)+\ka\cU(\rho_\ka)
 \leq Cr_\ka^2.
\end{equation}
Moreover, centered minimizers satisfy
\begin{equation}\label{zero-diffusion-compactness}
 d_\al(\rho_\ka,\cP_\Delta)\longrightarrow0
\end{equation}
uniformly over the choice of minimizer.  Their multipliers satisfy
\begin{equation}\label{multiplier-rate}
0 \leq \ell_\ka-\ell_\Delta \leq Cr_\ka^2, 
 \qquad \ell_\Delta:=2\cW(\rho_\Delta)=(1-u)w(1).
\end{equation}
\end{lemma}
\begin{proof}
Choose a non-negative $\vphi\in C_c^\infty(B_1)$ with
$\int\vphi=1$, and define the near-simplex measure
\[
 \nu_\ka(x):=u\sum_{i=0}^d r_\ka^{-d}
       \vphi\left(\frac{x-x_i}{r_\ka}\right).
\]
For small $\ka$, we have $r_\ka\leq1/2$, so the summands have disjoint supports.
Write $x=x_i+r_\ka z$, $y=x_j+r_\ka z'$ with $z,z'\in B_1$, so that
$x-y=e_{ij}+h$ with $|h|\le 2r_\ka\le1$.  For $i\ne j$, since
$|e_{ij}|=1$, the point $e_{ij}$ minimizes $W$, with $W(e_{ij})=w(1)$;
in particular $\grad W(e_{ij})=0$, and Taylor's formula gives
\[ 
0\le W(e_{ij}+h)-W(e_{ij})\le \tfrac12\Big(\max_{\overline{B_2}}|D^2W|\Big)|h|^2 \le Cr_\ka^2. 
\]
Integrating against $\nu_\ka\otimes\nu_\ka$, we find that the
cross-interaction excess of $\nu_\ka$ from $\rho_\Delta$ is at most $Cr_\ka^2$.
For $i=j$, $|x-y|\le 2r_\ka\le1$, so
$|W(x-y)|\le (2r_\ka)^\al/\al+(2r_\ka)^\bt/\bt\le Cr_\ka^\bt\le Cr_\ka^2$,
using $\al>\bt\ge2$.  Thus the self-interactions have
absolute value at most $Cr_\ka^2$.  Finally,
$\ka\cU(\nu_\ka) \le Cr_\ka^2$ by \eqref{cluster-scale-balance}.  Hence
\[
 \cE_\ka(\nu_\ka)\leq\cW(\rho_\Delta)+Cr_\ka^2.
\]
Minimality, $\cU\geq0$, and $\cW\geq\cW(\rho_\Delta)$ then prove
\eqref{energy-rate}. Since
\[
 0\le\ell_\ka-\ell_\Delta
 =\ka m\cU(\rho_\ka)+2\big(\cW(\rho_\ka)-\cW(\rho_\Delta)\big)
 \le 2m\big(\cE_\ka(\rho_\ka)-\cW(\rho_\Delta)\big),
\]
where both terms in the middle are non-negative, \eqref{energy-rate}
gives \eqref{multiplier-rate}.

Let $\ka_n\downarrow0$ and choose centered global minimizers.
Lemma \ref{existence} places them in a common compact ball, hence  a subsequence converges in $d_\al$ to a centered probability measure
$\rho_\infty$.  Continuity of the interaction energy on this ball and
\eqref{energy-rate} give $\cW(\rho_\infty)=\cW(\rho_\Delta)$, so
$\rho_\infty\in\cP_\Delta$ by \cite[Theorem 1.5]{DLM23}. 
If \eqref{zero-diffusion-compactness} failed uniformly, a sequence at
positive distance from $\cP_\Delta$ would have such a subsequence,
a contradiction.
\end{proof}

\begin{lemma}[Continuity of the interaction potential]\label{potential-convergence}
Let $R>0$, let $K\subset\R^d$ be compact, and let $V_\Delta:=W*\rho_\Delta$.
For every $\ep>0$ there is $\delta>0$ such that
\[
 \|W*\rho-V_\Delta\|_{C^2(K)}<\ep
\]
for every probability measure $\rho$ supported in $\overline{B_R}$ with
$d_\al(\rho,\rho_\Delta)<\delta$.
\end{lemma}
\begin{proof}
We may assume that $\rho_\Delta$ is supported in $\overline{B_R}$.  Let
$R':=R+\max_{x\in K}|x|$.  For $j=0,1,2$, the derivatives $D^jW$ are
bounded by some $M$ and uniformly continuous on $\overline{B_{R'}}$,
with a common modulus of continuity $\omega$.  Let $\pi$ be an optimal
coupling of $\rho$ and $\rho_\Delta$ for $d_\al$.  For $x\in K$ and $t>0$, split the integral below according to whether
$|y-y'|<t$ or $|y-y'|\geq t$.  On the first set we use the modulus
$\omega$, and on the second we use the bound $2M$ together with
Markov's inequality $\pi(|y-y'|\geq t)\leq t^{-\al}d_\al(\rho,\rho_\Delta)^\al$:
\[
 \begin{split}
 \big|D^j(W*\rho)(x)-D^jV_\Delta(x)\big|
 &\leq\int\big|D^jW(x-y)-D^jW(x-y')\big|\,d\pi(y,y')\\
 &\leq\omega(t)+2Mt^{-\al}d_\al(\rho,\rho_\Delta)^\al.
 \end{split}
\]
Choose $t$ with $\omega(t)<\ep/6$, and then $\delta$ with
$2Mt^{-\al}\delta^\al<\ep/6$.
\end{proof}

\subsection{Localization and the cluster scale}\label{subsection-scale}
We now use the energy bound \eqref{energy-rate} to determine the size
and location of the clusters.  The isolated, non-degenerate wells of
the simplex potential provide their local confinement.

\begin{lemma}[Non-degeneracy of the simplex wells]\label{simplex-wells}
Let \((\al,\bt)\in\cS\), let \(V_\Del:=W*\rho_\Del\), and let
$a_\perp,a_\parallel$ be as in \eqref{well-constants}.  Then
\[
 \argmin_{\R^d}V_\Del=\{x_0,\ldots,x_d\},
 \qquad
 \grad^2V_\Del(x_i)
 =u\Big(a_\perp\big(I_d-\hat x_i\otimes\hat x_i\big)
     +a_\parallel\,\hat x_i\otimes\hat x_i\Big).
\]
Moreover, $a_\parallel>0$, and $a_\perp>0$ if $d\geq2$.  In particular,
the constant
\[
 \sigma_\Delta:=\min_{0\leq i\leq d}
       \lambda_{\min}\bigl(\grad^2V_\Delta(x_i)\bigr)>0
\]
satisfies $\grad^2V_\Delta(x_i)\geq\sigma_\Delta I_d$ for
$i=0,\ldots,d$.
\end{lemma}

\begin{proof}
We first identify the minimum set. By \cite[Theorem 1.5]{DLM23}, $\rho_\Delta$ minimizes $\cW_{\al,\bt}$, hence the 
Euler--Lagrange inequality gives the inclusion $\{x_0,\ldots,x_d\}\subset\argmin V_\Delta$. To show the reverse inclusion, as $\al>\al_\Delta(\bt)\geq\bt$,  choose $\al'\in(\al_\Delta(\bt),\al)$.  By \cite[Theorem 1.5]{DLM23}, $\rho_\Delta$ also minimizes $\cW_{\al',\bt}$.  Since $\al>\al'$, the pair
$(\al,\bt)$ lies in the northeast region $A_{\al',\bt}$ of
\cite{DLM23}.  The northeast comparison principle
\cite[Theorem 1.2, (1.7)]{DLM23}, applied with $\nu=\rho_\Delta$ and
$(\al',\bt)$ in place of $(\al,\bt)$, therefore gives
\[
 \{x_0,\ldots,x_d\}=\supp\rho_\Delta=\argmin_{\R^d}\big(W_{\al,\bt}*\rho_\Delta\big)
 =\argmin_{\R^d}V_\Delta.
\]
For the reader's convenience we recall why the comparison principle
yields equality.  Write $\overline w_{a,b}(r):=(br^a-ar^b)/(a-b)$,
which is a positive multiple of $w_{a,b}$ plus a constant, so that
$\overline W_{a,b}*\rho_\Delta$ and $W_{a,b}*\rho_\Delta$ have the same
minimum points.  By \cite[Lemma 3.1]{DLM23} and the symmetry
$\overline w_{a,b}=\overline w_{b,a}$, for each fixed $r \ge 0$, the function
$a\mapsto\overline w_{a,\bt}(r)$ is nondecreasing for $a>\bt$, and
strictly increasing unless $r\in\{0,1\}$.  Hence
$\overline W_{\al',\bt}*\rho_\Delta\leq\overline W_{\al,\bt}*\rho_\Delta$,
with equality at $x$ if and only if $|x-x_j|\in\{0,1\}$ for every $j$.
A point with this property is a vertex: otherwise it would be at
distance one from all $d+1$ vertices, and $\R^d$ contains no $d+2$
points at mutual distance one.  Since the vertices minimize
$\overline W_{\al',\bt}*\rho_\Delta$ by the Euler--Lagrange equation for
$\rho_\Delta$, every non-vertex $x$ satisfies
$\overline W_{\al,\bt}*\rho_\Delta(x)>\overline W_{\al',\bt}*\rho_\Delta(x)
\geq\overline W_{\al',\bt}*\rho_\Delta(x_0)=\overline W_{\al,\bt}*\rho_\Delta(x_0)$.

We next compute the Hessian at each vertex.  For $i\neq j$,
\[
 \grad^2W(e_{ij})=(\al-\bt)e_{ij}\otimes e_{ij},
 \qquad
 \grad^2W(0)=-\ibt I_d.
\]
By \eqref{simplex-geometry},
\begin{equation*}
 \sum_{j\neq i}e_{ij}\otimes e_{ij}
 =(d+1)x_i\otimes x_i+\sum_jx_j\otimes x_j
 =\frac12I_d+\frac d2\,\hat x_i\otimes\hat x_i.
\end{equation*}
Consequently
\begin{equation}\label{simplex-well-hessian}
 \grad^2V_\Del(x_i) = u \Big(\sum_{j \neq i} \grad^2W(e_{ij}) + \grad^2W(0)\Big)
 =u\Big((\al-\bt)\sum_{j\neq i}e_{ij}\otimes e_{ij}-\ibt I_d\Big),
\end{equation}
which is the stated formula.  If $\bt>2$, then
$a_\parallel>a_\perp=(\al-\bt)/2>0$.  If $\bt=2$, then $(\al,2)\in\cS$
means $\al>4$ for $d\geq2$ and $\al>3$ for $d=1$; see
\cite[Theorem 1.5 and Remark 1.3]{DLM23}.  Thus $a_\perp=(\al-4)/2>0$
when $d\geq2$, and $a_\parallel=((\al-2)(d+1)-2)/2>0$ in both cases.
When $d=1$, the term containing $a_\perp$ vanishes.
\end{proof}

The following elementary rigidity statement places an almost regular
simplex near the reference simplex.

\begin{lemma}[Rigidity of almost regular simplices]\label{simplex-rigidity}
There are $\ep_R>0$ and $C_R$, depending only on $d$, with the following
property.  If $z_0,\ldots,z_d\in\R^d$ satisfy
$\big||z_i-z_j|-1\big|\leq\ep\leq\ep_R$ for all $i\neq j$, then there is a
rigid motion $T$ with $\max_i|Tz_i-x_i|\leq C_R\ep$.
\end{lemma}

\begin{proof}
Let $X$ and $Y$ be the $d\times d$ matrices with columns $x_i-x_0$ and
$z_i-z_0$, $1\leq i\leq d$. The polarization identity
$2a\cdot b=|a|^2+|b|^2-|a-b|^2$ expresses the entries of the Gram
matrices $X^TX$ and $Y^TY$ through the edge lengths, which are $1$ and
$1+O(\ep)$, respectively.  Hence $Y^TY=X^TX+E$ with $|E|\leq C\ep$.  Put
$M:=YX^{-1}$.  Then
\[
 M^TM=X^{-T}(X^TX+E)X^{-1}=I_d+X^{-T}EX^{-1}=I_d+O(\ep).
\]
For $\ep\leq\ep_R$ small, $M^TM$ is positive definite, so $M$ is
invertible and has the polar decomposition $M=QP$ with $Q\in O(d)$ and
$P=(M^TM)^{1/2}=I_d+O(\ep)$.  Consequently $|M-Q|=|Q(P-I_d)|\leq C\ep$
and $|Y-QX|\leq|M-Q|\,|X|\leq C\ep$, that is,
$|(z_i-z_0)-Q(x_i-x_0)|\leq C\ep$ for every $i$.  The rigid motion
$T(x):=Q^T(x-z_0)+x_0$ satisfies
$Tz_i-x_i=Q^T\big((z_i-z_0)-Q(x_i-x_0)\big)$, which proves the claim.
\end{proof}

We now show that every centered minimizer close to the simplex splits, after a rotation, into $d+1$ clusters of width $O(r_\ka)$ around the vertices, with densities of order $r_\ka^{-d}$ and masses within $O(r_\ka)$ of $u$.

\begin{lemma}[The scale of the clusters]\label{cluster-scale}
There exist \(\eta_0\in(0,1/8)\), \(\delta_{\mathrm{sc}}>0\), \(C>1\), and
\(\ka_{\mathrm{sc}}>0\) with the following property.  Let \(0<\ka<\ka_{\mathrm{sc}}\), and let
\(\rho_\ka\in\cP_0(\R^d)\cap\argmin\cE_\ka\) satisfy
\(d_\al(\rho_\ka,\rho_\Del)<\delta_{\mathrm{sc}}\).
Then $\supp\rho_\ka\subset\bigcup_iB_{\eta_0/2}(x_i)$, and there is
$Q\in O(d)$ such that the clusters $\rho_{\ka,i}$ of $Q_\#\rho_\ka$, with
masses $M_{\ka,i}:=|\rho_{\ka,i}|$, satisfy
\begin{equation}\label{cluster-scale-estimates}
 \frac{u}{2}\leq M_{\ka,i}\leq\frac{3u}{2},\qquad
 \supp\rho_{\ka,i}\subset B_{Cr_\ka}(x_i),\qquad
 \|\rho_{\ka,i}\|_{L^\infty}\leq Cr_\ka^{-d},
\end{equation}
and
\begin{equation}\label{mass-first-rate}
 \sum_{i=0}^{d}|M_{\ka,i}-u|^2\leq Cr_\ka^2.
\end{equation}
\end{lemma}

\begin{proof}
We localize the support in uniformly convex wells, use the mass
constraint to determine their heights and radii, and then use the
interaction excess to locate the wells relative to the reference simplex.

\proofstep{Step 1: Localization in the wells.}
Choose $\eta_0\in(0,1/8)$ so small that
$D^2V_\Delta\geq(3\sigma_\Delta/4)I_d$ on $\overline{B_{2\eta_0}(x_i)}$
for every $i$.  Lemma \ref{existence} places every centered minimizer in
a common ball $B_R$, and we may assume $R\geq1$, so that
$K:=\overline{B_R}$ contains these closed balls.  
By Lemma \ref{simplex-wells}, there is $g>0$ such that
$V_\Delta-\ell_\Delta\geq g$ on the compact set
$K\setminus\bigcup_iB_{\eta_0/2}(x_i)$, where
$\ell_\Delta=V_\Delta(x_i)$ is the minimum value of $V_\Delta$.  By
Lemma \ref{potential-convergence}, $\Phi_\ka$ is close to $V_\Delta$ in
$C^2(K)$ when $d_\al(\rho_\ka,\rho_\Delta)<\delta_{\mathrm{sc}}$ with
$\delta_{\mathrm{sc}}$ small, and $\ell_\ka\to\ell_\Delta$ uniformly over
minimizers by \eqref{multiplier-rate}.  First choose $\delta_{\mathrm{sc}}$,
and then $\ka_{\mathrm{sc}}$, so that
\[
 \Phi_\ka-\ell_\ka\geq\frac g2
 \quad\hbox{on }K\setminus\bigcup_iB_{\eta_0/2}(x_i),
 \qquad
 \|D^2\Phi_\ka-D^2V_\Delta\|\leq\frac{\sigma_\Delta}4
 \quad\hbox{on }K.
\]
Since $\supp\rho_\ka\subset K$, the first bound and
\eqref{rho-version} give
$\supp\rho_\ka\subset\bigcup_iB_{\eta_0/2}(x_i)$. 
The second bound gives, with $c_*:=\sigma_\Delta/2$ and some $C_*>c_*$
independent of the minimizer,
\begin{equation}\label{uniform-well-convexity}
 c_*I_d\leq\grad^2\Phi_\ka\leq C_*I_d
 \quad\hbox{on }B_{2\eta_0}(x_i),\qquad i=0,\ldots,d.
\end{equation}
By decreasing $\delta_{\mathrm{sc}}$ further, the masses of the clusters of
$\rho_\ka$ lie between $u/2$ and $3u/2$.  To see this, let $\pi$ be an
optimal coupling of $\rho_\ka$ and $\rho_\Delta$ for $d_\al$, and write
$B_i:=B_{\eta_0/2}(x_i)$.  Since $\rho_\ka$ is concentrated on
$\bigcup_jB_j$, we have
\[
 M_{\ka,i}-u=\pi\big(B_i\times\{x_j:j\neq i\}\big)
 -\pi\Big(\bigcup_{j\neq i}B_j\times\{x_i\}\Big).
\]
Both sets on the right consist of pairs at distance at least
$1-\eta_0/2\geq1/2$, so Markov's inequality gives
$|M_{\ka,i}-u|\leq2^\al d_\al(\rho_\ka,\rho_\Delta)^\al$.  It suffices
that $2^\al\delta_{\mathrm{sc}}^\al\leq u/2$.

\proofstep{Step 2: Heights and radii of the clusters.}
The \(C^1\)-closeness to \(V_\Del\), together with
\eqref{uniform-well-convexity}, gives a unique critical point
\(z_{\ka,i}\in B_{\eta_0/2}(x_i)\), the minimum of $\Phi_\ka$ in
$B_{2\eta_0}(x_i)$.  Indeed, since $\grad V_\Delta(x_i)=0$ and
$D^2V_\Delta\geq(3\sigma_\Delta/4)I_d$ on $B_{2\eta_0}(x_i)$, every
$x$ with $|x-x_i|=\eta_0/2$ satisfies
\[
 (x-x_i)\cdot\grad\Phi_\ka(x)
 \geq\frac{3\sigma_\Delta}{4}\frac{\eta_0^2}{4}
   -\frac{\eta_0}{2}\,\|\grad\Phi_\ka-\grad V_\Delta\|_{L^\infty(K)}>0
\]
if the $C^1$ distance is below $3\sigma_\Delta\eta_0/16$.  Thus the
outward radial derivative of the strictly convex function $\Phi_\ka$ is
positive on $\partial B_{\eta_0/2}(x_i)$, its minimum over the closed
ball $\overline{B_{2\eta_0}(x_i)}$ lies in $B_{\eta_0/2}(x_i)$, and
strict convexity gives uniqueness of the critical point.  Put
\[
 h_{\ka,i}:=\ell_\ka-\Phi_\ka(z_{\ka,i}).
\]
By \eqref{rho-version},
\begin{equation}\label{density-from-potential}
 \rho_\ka(x)=
 \left[\frac{m-1}{m\ka}\big(\ell_\ka-\Phi_\ka(x)\big)\right]_+^{1/(m-1)}.
\end{equation}
Since the \(i\)-th cluster has positive mass, \(h_{\ka,i}>0\).
For \(x\in B_{2\eta_0}(x_i)\),
\eqref{uniform-well-convexity} implies
\[
 \frac{c_*}{2}|x-z_{\ka,i}|^2
 \leq \Phi_\ka(x)-\Phi_\ka(z_{\ka,i})
 \leq\frac{C_*}{2}|x-z_{\ka,i}|^2.
\]

First, \(h_{\ka,i}\to0\) uniformly over the minimizers in
question: Otherwise there are $h_0>0$ and minimizers with
$\ka\downarrow0$ and $h_{\ka,i}\geq h_0$.  The upper quadratic bound
gives $\ell_\ka-\Phi_\ka\geq h_0/2$ on $B_{(h_0/C_*)^{1/2}}(z_{\ka,i})$,
which lies in $B_{2\eta_0}(x_i)$ if $h_0$ is small.  By
\eqref{density-from-potential}, the mass of this ball is at least
$c\,(h_0/C_*)^{d/2}\big(h_0/\ka\big)^{1/(m-1)}\to\infty$, contradicting
$M_{\ka,i}\leq3u/2$.

By the two quadratic bounds, on $B_{2\eta_0}(x_i)$ we have
\[
 h_{\ka,i}-\frac{C_*}{2}|x-z_{\ka,i}|^2
 \leq\ell_\ka-\Phi_\ka(x)
 \leq h_{\ka,i}-\frac{c_*}{2}|x-z_{\ka,i}|^2 .
\]
Hence $\ell_\ka-\Phi_\ka\geq h_{\ka,i}/2$ on
$B_{(h_{\ka,i}/C_*)^{1/2}}(z_{\ka,i})$, while
$[\ell_\ka-\Phi_\ka]_+\leq h_{\ka,i}$ on $B_{2\eta_0}(x_i)$ and
$\{\ell_\ka>\Phi_\ka\}\cap B_{2\eta_0}(x_i)\subset
B_{(2h_{\ka,i}/c_*)^{1/2}}(z_{\ka,i})$.  Since $h_{\ka,i}\to0$, both balls
lie in $B_{2\eta_0}(x_i)$ for small $\ka$.  Since
$\supp\rho_\ka\subset\bigcup_jB_{\eta_0/2}(x_j)$, the mass of the $i$-th
cluster is $M_{\ka,i}=\rho_\ka(B_{2\eta_0}(x_i))$.  Integrating
\eqref{density-from-potential} over the two balls therefore gives
\begin{equation}\label{height-mass-comparison}
 c\ka^{-1/(m-1)}h_{\ka,i}^{\frac d2+\frac{1}{m-1}}
 \leq M_{\ka,i}\leq
 C\ka^{-1/(m-1)}h_{\ka,i}^{\frac d2+\frac{1}{m-1}},
\end{equation}
with constants depending only on $c_*$, $C_*$, $d$, and $m$.  Since the
cluster masses lie between $u/2$ and $3u/2$, \eqref{height-mass-comparison} gives
\begin{equation}\label{height-scale}
 cr_\ka^2\leq h_{\ka,i}\leq Cr_\ka^2.
\end{equation}
Equations \eqref{density-from-potential}--\eqref{height-scale} now imply that the $i$-th cluster of $\rho_\ka$ is supported in $B_{Cr_\ka}(z_{\ka,i})$ and is bounded by $\big[\frac{m-1}{m\ka}h_{\ka,i}\big]^{1/(m-1)}\leq Cr_\ka^{-d}$, since $r_\ka^2/\ka=r_\ka^{-d(m-1)}$ by \eqref{cluster-scale-balance}.

\proofstep{Step 3: Masses and positions of the clusters.}
The radius estimate is currently centered at $z_{\ka,i}$, not at the
reference vertex $x_i$.  We now extract the information needed to place
these centers.  Let $\rho_{\ka,i}$ denote the clusters of $\rho_\ka$,
with masses $M_{\ka,i}$, and put
\[
 I_{\ka,ij}:=\iint W(x-y)\,d\rho_{\ka,i}(x)\,d\rho_{\ka,j}(y).
\]
The preceding support estimate gives \(|I_{\ka,ii}|\leq Cr_\ka^\bt\leq Cr_\ka^2\), since $\bt\geq2$ and $r_\ka\leq1$. For \(i\neq j\), define the cross excess
\[
 X_{\ka,ij}:=I_{\ka,ij}-w(1)M_{\ka,i}M_{\ka,j}\geq0,
\]
which is non-negative because \(w(1)=\min W\).  Since $\rho_\ka=\sum_i\rho_{\ka,i}$ and $W$ is even,
$\cW(\rho_\ka)=\frac12\sum_iI_{\ka,ii}+\sum_{i<j}I_{\ka,ij}$.  Using
$\sum_iM_{\ka,i}=1$, we have
$\sum_{i<j}M_{\ka,i}M_{\ka,j}=\frac12\big(1-\sum_iM_{\ka,i}^2\big)$
and $\sum_iM_{\ka,i}^2-u=\sum_i(M_{\ka,i}-u)^2$. With
$\cW(\rho_\Del)=\frac{w(1)}2(1-u)$, this gives
\begin{equation*}
 \cW(\rho_\ka)-\cW(\rho_\Del)
 =\frac{c_W}{2}\sum_i(M_{\ka,i}-u)^2
   +\sum_{i<j}X_{\ka,ij}
   +\frac12\sum_i I_{\ka,ii}
\end{equation*}
where $c_W = -w(1) >0$. Combining this identity with \eqref{energy-rate}, we get
\eqref{mass-first-rate} and
\begin{equation}\label{cross-excess-rate}
 0\leq X_{\ka,ij}\leq Cr_\ka^2.
\end{equation}

It remains to find $Q\in O(d)$ with \(\max_i|Qz_{\ka,i}-x_i|\leq Cr_\ka\). Indeed, since $\supp\rho_{\ka,i}\subset B_{Cr_\ka}(z_{\ka,i})$, such a $Q$ gives $\supp(Q_\#\rho_{\ka,i})\subset B_{2Cr_\ka}(x_i)\subset B_{1/2}(x_i)$ for small $\ka$, so $Q_\#\rho_{\ka,i}$ is the $i$-th cluster of $Q_\#\rho_\ka$, with the same mass, and \eqref{cluster-scale-estimates} and \eqref{mass-first-rate} hold for $Q_\#\rho_\ka$.

Since \(w\) has a non-degenerate minimum at \(s=1\), there is
\(c>0\) such that, whenever \(s\) is close to one,
$w(s)-w(1)\geq c(s-1)^2$. For \(x\in\supp(\rho_{\ka,i})\) and
\(y\in\supp(\rho_{\ka,j})\), the preceding support estimate gives
\(
 ||x-y|-|z_{\ka,i}-z_{\ka,j}||\leq Cr_\ka
\).  Hence
\[
 W(x-y)-w(1)
 \geq c\big(|z_{\ka,i}-z_{\ka,j}|-1\big)^2-Cr_\ka^2.
\]
After integration by $\rho_{\ka,i}\otimes\rho_{\ka,j}$, \eqref{cross-excess-rate} and the mass bound $\frac12 u \le M_{\ka,i} \le \frac32 u$ give
\[ 
\big||z_{\ka,i}-z_{\ka,j}|-1\big|\leq Cr_\ka,
 \qquad i\neq j.
\]
For small
$\ka$, Lemma~\ref{simplex-rigidity}, applied with $\ep=Cr_\ka$, provides
a rigid motion \(T_\ka(x)=Q_\ka x+a_\ka\) such that
$\max_i|T_\ka(z_{\ka,i})-x_i|\leq Cr_\ka$. 
Since \(\rho_\ka\) is centered and $\supp \rho_{\ka,i} \subset B_{Cr_\ka}(z_{\ka,i})$,
\[
 0=\sum_i\int x\,d\rho_{\ka,i}(x)
   =\sum_iM_{\ka,i}z_{\ka,i}+O(r_\ka).
\]
Since $\sum_iM_{\ka,i}=1$, this gives
\[
 \sum_iM_{\ka,i}T_\ka(z_{\ka,i})
 =Q_\ka\sum_iM_{\ka,i}z_{\ka,i}+a_\ka=a_\ka+O(r_\ka).
\]
On the other hand, by the estimate for $T_\ka(z_{\ka,i})$, \eqref{mass-first-rate}, and $\sum_ix_i=0$,
\[
 \sum_iM_{\ka,i}T_\ka(z_{\ka,i})=\sum_iM_{\ka,i}x_i+O(r_\ka)
 =\sum_i(M_{\ka,i}-u)x_i+O(r_\ka)=O(r_\ka).
\]
Hence \(|a_\ka|\leq Cr_\ka\), and therefore 
\(\max_i|Q_\ka z_{\ka,i}-x_i|\leq Cr_\ka\).\end{proof}

\subsection{Alignment and normalized minimizers}\label{subsection-alignment}
The scale estimates leave the rigid-motion freedom unresolved.  We fix
it by the first-order conditions for the best rigid fit of the clusters
to the labelled vertices.

\begin{definition}[Alignment]\label{alignment-definition}
Let $\rho_0,\ldots,\rho_d$ be finite non-negative measures with finite
second moments.  The labelled family $(\rho_i)$, or the measure
$\rho=\sum_i\rho_i$ with these clusters, is \emph{aligned with}
$\Delta^d$ if
\begin{equation}\label{alignment-conditions}
 \sum_{i=0}^d\int(x-x_i)\,d\rho_i(x)=0,
 \qquad
 \sum_{i=0}^d\int(x-x_i)\cdot Sx_i\,d\rho_i(x)=0
\end{equation}
for every skew-symmetric matrix $S$.
\end{definition}

The conditions \eqref{alignment-conditions} say that the identity is a
critical point of
\[
 J(T):=\sum_{i=0}^d\int|Tx-x_i|^2\,d\rho_i(x)
\]
on the group of rigid motions.  Indeed, for a skew-symmetric $S$ and $b\in\R^d$, the maps $R_t(x):=e^{tS}x+tb$ form a $C^1$ curve of rigid motions with $R_0=\mathrm{id}$ and velocity $x\mapsto Sx+b$ at $t=0$, because $(e^{tS})^Te^{tS}=e^{-tS}e^{tS}=I_d$. Conversely, if $t\mapsto Q_tx+a_t$ is any $C^1$ curve of rigid motions with $Q_0=I_d$ and $a_0=0$, differentiating $Q_t^TQ_t=I_d$ at $t=0$ shows that $\dot Q_0$ is skew-symmetric, so its velocity $x\mapsto\dot Q_0x+\dot a_0$ is again of this form. Since the first variation of $J$ along a curve depends only on its velocity at $t=0$, it suffices to consider the curves $R_t$. 
Using $(x-x_i)\cdot S(x-x_i)=0$,
differentiation at $t=0$ gives
\begin{equation*}
 \frac{d}{dt}\Big|_{t=0}J(R_t)
 =2\,b\cdot\sum_i\int(x-x_i)\,d\rho_i(x)
  +2\sum_i\int(x-x_i)\cdot Sx_i\,d\rho_i(x).
\end{equation*}
This vanishes for all $(S,b)$ if and only if
\eqref{alignment-conditions} holds. Note that the conditions involve only the
masses $M_i:=|\rho_i|$ and the first moments
$\mathbf{m}_i:=\int x\,d\rho_i$, and since $v\cdot Sw=\langle S,v\otimes w\rangle$,
the second condition says exactly that
$\sum_i(\mathbf{m}_i-M_ix_i)\otimes x_i$ is symmetric.

\begin{lemma}[Fixing translation and orientation]\label{alignment-lemma}
Let $\rho=\sum_i\rho_i$ be a probability measure with finite second
moment, with labelled clusters $\rho_i$, masses $M_i$, and first
moments $\mathbf{m}_i$, and put
\[
 \bar x:=\sum_iM_ix_i,\qquad \bar{\mathbf m}:=\sum_i\mathbf{m}_i,\qquad
 K:=\sum_i\mathbf{m}_i\otimes(x_i-\bar x).
\]
\begin{enumerate}
\item[\textup{(i)}] A rigid motion $x\mapsto Q'x+a'$ aligns the clusters
if and only if $a'=\bar x-Q'\bar{\mathbf m}$ and $Q'K$ is symmetric.  If
$K$ is invertible, then $Q:=(K^TK)^{1/2}K^{-1}$ is orthogonal, $QK$ is
symmetric and positive definite, and $T(x):=Qx+a$ with
$a:=\bar x-Q\bar{\mathbf m}$ aligns the clusters.  Moreover, $T$ is the
unique global minimizer of $J$ over all rigid motions, and if
$\|K-\frac u2I_d\|\leq\frac u8$, then every aligning rigid motion
$x\mapsto Q'x+a'$ with $\|Q'-I_d\|\leq\frac18$ coincides with $T$.
\item[\textup{(ii)}] For every $C>0$ there are $r_0,C'>0$ such that, if
$0<r\leq r_0$, $\supp\rho_i\subset B_{Cr}(x_i)$, and
$|M_i-u|\leq Cr$ for all $i$, then $K$ is invertible and
$\|Q-I_d\|+|a|\leq C'r$.
\item[\textup{(iii)}] If $G$ is a symmetry of $\Delta^d$ with
$Gx_i=x_{\pi(i)}$ and $(\rho_i)$ is aligned, then the relabelled family
$(\widetilde\rho_j)$, $\widetilde\rho_{\pi(i)}:=G_\#\rho_i$, is aligned.
\end{enumerate}
\end{lemma}

\begin{proof}
(i) For a rigid motion $T'(x)=Q'x+a'$, expand the square and use
$\sum_iM_i=1$, $|Q'x|=|x|$, and
$(Q'\mathbf{m}_i)\cdot(x_i-\bar x)=\tr\big(Q'\,\mathbf{m}_i\otimes(x_i-\bar x)\big)$.
Completing the square in $a'$ gives
\begin{equation}\label{J-reduced}
 J(T')=c_0-2\tr(Q'K)+\big|a'-\bar x+Q'\bar{\mathbf m}\big|^2,
\end{equation}
where
$c_0:=\sum_i\int|x|^2\,d\rho_i+\sum_iM_i|x_i|^2-|\bar{\mathbf m}|^2-|\bar x|^2$
does not depend on $T'$.

The clusters $T'_\#\rho_i$ are aligned if and only if the identity is a
critical point of $R\mapsto\sum_i\int|Rx-x_i|^2\,d(T'_\#\rho_i)=J(R\circ T')$.
Along $R_t(x)=e^{tS}x+tb$, the motion $R_t\circ T'$ has linear part
$e^{tS}Q'$ and translation $e^{tS}a'+tb$.  Put
$\delta:=a'-\bar x+Q'\bar{\mathbf m}$.  By \eqref{J-reduced},
\[
 \frac{d}{dt}\Big|_{t=0}J(R_t\circ T')
 =-2\tr(SQ'K)+2\,\delta\cdot\big(S(a'+Q'\bar{\mathbf m})+b\big).
\]
Taking $S=0$ and $b$ arbitrary shows that alignment forces $\delta=0$.
When $\delta=0$, the derivative is $-2\tr(SQ'K)=2\langle S,Q'K\rangle$,
which vanishes for every skew-symmetric $S$ if and only if $Q'K$ is
symmetric.  This proves the characterization.

If $K$ is invertible, then $P:=(K^TK)^{1/2}$ is symmetric and positive
definite, $Q^TQ=K^{-T}P^2K^{-1}=I_d$, and $QK=P$.  The characterization
shows that $T$ aligns the clusters.  To see that $T$ minimizes $J$, note
by \eqref{J-reduced} that it suffices to maximize $\tr(Q'K)$ over
$Q'\in O(d)$.  Write $O:=Q'Q^T\in O(d)$, so that $\tr(Q'K)=\tr(OP)$, and
let $(e_j)$ be an orthonormal eigenbasis of $P$ with eigenvalues
$\lambda_j>0$.  Then
\[
 \tr(OP)=\sum_j\lambda_j\,e_j\cdot Oe_j\leq\sum_j\lambda_j=\tr(QK),
\]
with equality only if $Oe_j=e_j$ for all $j$, that is, $Q'=Q$.  Thus
$T$ uniquely minimizes $J$.

Finally, let $x\mapsto Q'x+a'$ align the clusters, with
$\|Q'-I_d\|\leq\frac18$, and suppose $\|K-\frac u2I_d\|\leq\frac u8$.
Then $\|K\|\leq\frac{5u}8$ and
$\|Q'K-\frac u2I_d \|
 \leq\|Q'-I_d\|\,\|K\|+ \|K-\frac u2I_d \|<\frac u4$. 
So the symmetric matrix $Q'K$ is positive definite.  Hence
$K=Q'^T(Q'K)=Q^TP$ are two polar decompositions of $K$, and uniqueness
of the polar decomposition gives $Q'=Q$.  Then $a'=a$ by the
characterization.

(ii) Under these assumptions, $|\mathbf{m}_i-M_ix_i|\leq Cr$,
$\bar x=\sum_i(M_i-u)x_i=O(r)$, and
$\bar{\mathbf m}=\sum_i(\mathbf{m}_i-M_ix_i)+\bar x=O(r)$.  The identity $ \sum_{j=0}^dx_j\otimes x_j=\frac12I_d$ in
\eqref{simplex-geometry} then gives
\[
 K=\sum_i(\mathbf{m}_i-M_ix_i)\otimes x_i+\sum_iM_ix_i\otimes x_i
   -\bar{\mathbf m}\otimes\bar x
 =\frac u2I_d+O(r).
\]
The map $K\mapsto(K^TK)^{1/2}K^{-1}$ is smooth near $\frac u2I_d$ and
equals $I_d$ there.  Hence $\|Q-I_d\|\leq C'r$, and then
$|a|\leq|\bar x|+|\bar{\mathbf m}|\leq C'r$.

(iii) The measures $G_\#\rho_i$ have finite second moments, and
\[
 \sum_i\int(y-x_{\pi(i)})\,d\widetilde\rho_{\pi(i)}(y)
 =G\sum_i\int(x-x_i)\,d\rho_i(x)=0.
\]
For a skew-symmetric $S$, the rotational condition for the relabelled
family becomes
$\sum_i\int(x-x_i)\cdot(G^TSG)x_i\,d\rho_i(x)=0$, which holds because
$G^TSG$ is again skew-symmetric.
\end{proof}

Alignment alone does not determine the rigid motion: by
Lemma~\ref{alignment-lemma}\textup{(i)}, every $Q'\in O(d)$ with $Q'K$
symmetric gives an aligning motion, and only $Q'=Q$ makes $Q'K$ positive
definite.  This is why the uniqueness statement in
Lemma~\ref{alignment-lemma}\textup{(i)} requires $Q'$ to be close to
$I_d$.

\begin{definition}[Normalized minimizers]\label{normalized-definition}
Let $C$ be the constant of Lemma \ref{cluster-scale}, let $C'$ be given
by Lemma \ref{alignment-lemma}\textup{(ii)} for this $C$, and put
$L:=C+2C'+C^{1/2}$.  For $\ka>0$ with $Lr_\ka<1/2$, let $\cN_\ka$ be the
set of global minimizers $\rho$ of $\cE_\ka$ such that
$\supp\rho\subset\bigcup_iB_{Lr_\ka}(x_i)$, the clusters of $\rho$ are
aligned with $\Delta^d$, and
\begin{equation}\label{normalized-estimates}
 \big|\,|\rho_i|-u\,\big|\leq Lr_\ka,
 \qquad
 \|\rho_i\|_{L^\infty}\leq Lr_\ka^{-d},
 \qquad i=0,\ldots,d.
\end{equation}
\end{definition}

\begin{lemma}[Normalization]\label{normalization-lemma} There is $\ka_3>0$ such that, for $0<\ka<\ka_3$, every global minimizer of $\cE_\ka$ is the image of an element of $\cN_\ka$ under a rigid motion; in particular, $\cN_\ka\neq\emptyset$. Moreover, $G_\#\rho\in\cN_\ka$ for every $\rho\in\cN_\ka$ and every symmetry $G$ of $\Delta^d$. \end{lemma}

\begin{proof}
Choose small $\ka_3\in(0,\ka_{\mathrm{sc}}]$ so that, for $\ka \in (0,\ka_3)$,  every centered global minimizer
satisfies $d_\al(\rho,\cP_\Delta)<\delta_{\mathrm{sc}}$ by Lemma \ref{compactsupport}, $r_\ka\leq r_0$
in Lemma \ref{alignment-lemma}\textup{(ii)}, and $Lr_\ka<1/2$.  Let
$\rho$ be a centered global minimizer with $0<\ka<\ka_3$. Since $\cP_\Delta$ is the $O(d)$-orbit of
$\rho_\Delta$, an orthogonal map then brings $\rho$ within
$d_\al$-distance $\delta_{\mathrm{sc}}$ of $\rho_\Delta$.
Lemma \ref{cluster-scale} provides a further orthogonal map after which
\eqref{cluster-scale-estimates} and \eqref{mass-first-rate} hold; in
particular, the clusters are compactly supported and
$|M_{\ka,i}-u|\leq C^{1/2}r_\ka$.  Lemma \ref{alignment-lemma}\textup{(ii)}
then gives a rigid motion $T(x)=Qx+a$ with $\|Q-I_d\|+|a|\leq C'r_\ka$
that aligns the clusters.  It maps $B_{Cr_\ka}(x_i)$ into
$B_{(C+2C')r_\ka}(x_i)$, and it preserves masses and $L^\infty$ bounds.
Hence $T_\#\rho\in\cN_\ka$.

Now let $\rho\in\cN_\ka$, and let $G$ be a symmetry of $\Delta^d$. Then 
$G_\#\rho$ is a global minimizer and the map $G$ permutes the balls
$B_{Lr_\ka}(x_i)$ and the clusters, so \eqref{normalized-estimates}
holds for $G_\#\rho$, which is aligned by
Lemma \ref{alignment-lemma}\textup{(iii)}.  Thus $G_\#\rho\in\cN_\ka$.
\end{proof}

Alignment is not the same as centering when cluster masses are unequal:
its first condition says $\int x\,d\rho=\sum_i|\rho_i|x_i$. Because of this, the
arguments below use the support and scale estimates in
\eqref{normalized-estimates}, not membership in $\cP_0(\R^d)$.  Once equal masses are proved, normalized minimizers are centered.

\subsection{Uniqueness of normalized minimizers}\label{subsection-uniqueness}
For a measure $\rho$ with clusters $\rho_i$ and for $r>0$, we call
\[
 f_i(z):=r^d\rho_i(x_i+rz)
\]
its rescaled clusters. We also define the rescaled kernels
\begin{align}\label{scaled-kernels}
 K_{ii}^r(z,z'):=\frac{W(r(z-z'))}{r^2},
 \quad 
 K_{ij}^r(z,z'):=\frac{W(e_{ij}+r(z-z'))-w(1)}{r^2} \ \ \,
 (i\neq j).
\end{align}
Taylor's formula, $W(0)=0$, $\grad W(0)=\grad W(e_{ij})=0$, and
continuity of $D^2W$ show that these kernels converge in $C^2$ on every
bounded set as \(r\downarrow0\) to
\begin{equation}\label{limiting-kernels}
 K_{ii}^0(z,z')=-\frac12\ibt|z-z'|^2,
 \quad
 K_{ij}^0(z,z')=\frac{\al-\bt}{2}
 \big(e_{ij}\cdot(z-z')\big)^2\quad(i\neq j).
\end{equation}

If a perturbation $(\delta_i)$ has zero mass in each cluster, then its limiting interaction, computed with the kernels \eqref{limiting-kernels}, depends only on the first moments $p_i:=\int z\,\delta_i(z)\,dz$: it equals the quadratic form $\cB(p)$ introduced below, as the proof of Lemma~\ref{local-uniqueness} shows. The next lemma shows $\cB(p)\geq0$ whenever $p$ satisfies the alignment conditions \eqref{moment-alignment}.

\begin{lemma}[Nonnegativity on aligned first moments]\label{dipole-positivity}
For $p=(p_0,\ldots,p_d)\in(\R^d)^{d+1}$ define
$\displaystyle \cB(p):=\ibt\sum_i|p_i|^2
       -2(\al-\bt)\sum_{i<j}(e_{ij}\cdot p_i)(e_{ij}\cdot p_j)$ and 
$\displaystyle M_p:=\sum_ip_i\otimes x_i$.
If
\begin{equation}\label{moment-alignment}
 \sum_ip_i=0
 \qquad\hbox{and}\qquad
 \sum_ip_i\cdot Sx_i=0\quad\hbox{for every skew-symmetric }S,
\end{equation}
then $M_p$ is symmetric and
\begin{equation}\label{dipole-identity}
 \cB(p)=\ibt\sum_i|p_i|^2+(\al-\bt)\Big((\tr M_p)^2+|M_p|^2\Big)\geq0,
\end{equation}
with equality only for $p=0$.
\end{lemma}

\begin{proof}
Since $p_i\cdot Sx_i=\langle S,p_i\otimes x_i\rangle$, the second
condition in \eqref{moment-alignment} says that $M_p$ is orthogonal to
all skew-symmetric matrices, that is, $M_p$ is symmetric.  Put
$s_{ij}:=x_i\cdot p_j$.  The first condition gives $\sum_js_{ij}=0$ for
every $i$.  Since $e_{ij}\cdot p_i=s_{ii}-s_{ji}$ and
$e_{ij}\cdot p_j=s_{ij}-s_{jj}$, and since the summand below is
symmetric in $i$ and $j$,
\[
 2\sum_{i<j}(e_{ij}\cdot p_i)(e_{ij}\cdot p_j)
 =\sum_{i\neq j}\big(s_{ii}s_{ij}-s_{ii}s_{jj}-s_{ji}s_{ij}+s_{ji}s_{jj}\big).
\]
Using $\sum_{j\neq i}s_{ij}=-s_{ii}$ in the first and last terms, we find
\[
 2\sum_{i<j}(e_{ij}\cdot p_i)(e_{ij}\cdot p_j)
 =-\Big(\sum_is_{ii}\Big)^2-\sum_{i,j}s_{ij}s_{ji}
 =-(\tr M_p)^2-\tr(M_p^2).
\]
As $M_p$ is symmetric, $\tr(M_p^2)=|M_p|^2$, which proves
\eqref{dipole-identity}.  If $\cB(p)=0$, then $M_p=0$ because $\al>\bt$,
so $\sum_ip_i\otimes(1,x_i)=0$ in $\R^d\otimes\R^{d+1}$.  The vectors
$(1,x_i)\in\R^{d+1}$ are linearly independent, since the vertices are
affinely independent, and therefore $p=0$.
\end{proof}

\begin{lemma}[Convex wells and inverse diffusion weights]\label{inverse-weight-lemma}
There are $\ka_2\in(0,\ka_3]$ and $C_w > 0$ such that, for
$0<\ka<\ka_2$ and every $\rho\in\cN_\ka$,
\begin{equation}\label{normalized-convexity}
 \grad^2(W*\rho)\geq c_*I_d\quad\hbox{on }B_{2\eta_0}(x_i),
 \qquad i=0,\ldots,d,
\end{equation}
where $\eta_0$ is as in Lemma \ref{cluster-scale} and
$c_*=\sigma_\Delta/2$ is the constant in \eqref{uniform-well-convexity}.
The rescaled clusters of $\rho$ with $r=r_\ka$ also satisfy
\begin{equation}\label{inverse-weight-bound}
 \int_{\{f_i>0\}}f_i^{2-m}\,dz\leq C_w,\qquad i=0,\ldots,d.
\end{equation}
\end{lemma}

\begin{proof}
Let $\rho\in\cN_\ka$, $r:=r_\ka$, $\Phi:=W*\rho$, and
$\ell:=\ka m\cU(\rho)+2\cW(\rho)$.  We take $\ka_2$ so small that
$Lr_\ka\leq\min\{\eta_0,u/2\}$ for $\ka<\ka_2$; in particular
$\supp\rho\subset B_1$, and every cluster has mass at least $u/2$.
To obtain an upper bound for $ d_\al(\rho,\rho_\Delta)$, we may move each cluster to its vertex, and then move the excess masses along edges of length one. 
This gives
\[
 d_\al(\rho,\rho_\Delta)^\al \leq (Lr)^\al + \sum_i\big||\rho_i|-u\big|(1+Lr)^\al
 \leq (Lr)^\al + \sum_i (Lr)(1+Lr)^\al \le Cr.
\]
Recall from Step 1 of the proof of Lemma \ref{cluster-scale} that
$D^2V_\Delta\geq(3\sigma_\Delta/4)I_d$ on $\overline{B_{2\eta_0}(x_i)}$.
Lemma \ref{potential-convergence}, applied with $R=1$ and
$K:=\bigcup_i\overline{B_{2\eta_0}(x_i)}$, therefore gives
\eqref{normalized-convexity} with the same $c_*=\sigma_\Delta/2$ for small enough $\ka_2$.

Set $D_r:=B_{2\eta_0/r}$. For $z\in D_r$, we have $f_i(z)=r^d\rho(x_i+rz)$, so \eqref{rho-version} gives 
\begin{equation}\label{eq123}
 f_i(z)=\Big[\frac{m-1}{m}\,\frac{r^{d(m-1)}}{\ka} \big(\ell-\Phi(x_i+rz)\big)\Big]_+^{1/(m-1)} =\Big[\frac{m-1}{m}P_i(z)\Big]_+^{1/(m-1)},
\end{equation}
where $P_i(z):= (\ell-\Phi(x_i+rz))/r^2$ and $r^{d(m-1)+2}=\ka$ from \eqref{cluster-scale-balance}. Since $D_z^2[\Phi(x_i+rz)]=r^2D^2\Phi(x_i+rz)$, \eqref{normalized-convexity} gives $D^2P_i\leq-c_*I_d$ on $D_r$. 
All level and superlevel sets below are taken inside $D_r$.
The rescaled clusters satisfy $\supp f_i\subset\overline{B_L}$,
$f_i\leq L$, and $\int f_i\geq u/2$. 
Since $P_i$ is nonpositive on $\overline{D_r}\setminus B_L$,
and positive somewhere in $\overline{B_L}$ because $\int f_i\geq u/2$, its maximum over
$\overline{D_r}$ is positive and attained in $\overline{B_L}$. 
Hence
$h_i:=\max_{\overline{D_r}}P_i$ satisfies
\[
 \frac u2\leq\int f_i\leq|B_L|\Big(\frac{m-1}{m}h_i\Big)^{1/(m-1)},
\]
so that $h_i\geq c_h$ for some $c_h>0$, independent of $\ka$.

If $1<m\leq2$, then $f_i^{2-m}\leq L^{2-m}$, and
\eqref{inverse-weight-bound} is immediate.  Let $m>2$, and put
$\theta:=(m-2)/(m-1)\in(0,1)$, so that
$f_i^{2-m}=((m-1)/m)^{-\theta}P_i^{-\theta}$ on $\{f_i>0\}$.  Fix a
maximum point $z^*$ of $P_i$ and $0<t\leq c_h/2$.  By concavity, the
superlevel set $\{P_i\geq t\}$ is convex.  It is contained in
$\overline{B_L}$, because $f_i>0$ there, and its boundary is the level
set $\{P_i=t\}$.  Concavity also gives
$P_i(z^*)\leq P_i(z)+\grad P_i(z)\cdot(z^*-z)$, hence
\[
 |\grad P_i|\geq\frac{h_i-t}{2L}\geq\frac{c_h}{4L}
 \quad\hbox{on }\{P_i=t\}.
\]
The perimeter of a convex set does not exceed that of a convex set
containing it, so
$\mathcal H^{d-1}(\{P_i=t\})\leq\mathcal H^{d-1}(\partial B_L)$; for
$d=1$, $\mathcal H^0$ is the counting measure.  The coarea formula $\int g\,dx=\int_{-\infty}^{\infty}(\int_{\{f=t\}}\frac{g}{|\nabla f|}\,d\mathcal H^{d-1})dt$ with $g = P_i^{-\theta}$ and $f = P_i$ yields
\[
 \begin{split}
 \int_{D_r\cap\{0<P_i<c_h/2\}}P_i^{-\theta}\,dz
 &=\int_0^{c_h/2}t^{-\theta}\int_{D_r\cap\{P_i=t\}}
   \frac{d\mathcal H^{d-1}}{|\grad P_i|}\,dt\\
 &\leq\frac{4L}{c_h}\,\mathcal H^{d-1}(\partial B_L)
   \int_0^{c_h/2}t^{-\theta}\,dt<\infty,
 \end{split}
\]
while $\int_{D_r\cap\{P_i\geq c_h/2\}}P_i^{-\theta}\,dz\leq(c_h/2)^{-\theta}|B_L|$.
This proves \eqref{inverse-weight-bound}.
\end{proof}

The bound \eqref{inverse-weight-bound} says that the inverse of the diffusion weight $U''(f_i) = mf_i^{m-2}$ is integrable on the positivity set of each cluster, uniformly in $\ka$. The proof of the following lemma needs this when $m>2$, where $U''$ degenerates at the free boundary.

\begin{lemma}[Uniqueness of normalized minimizers]\label{local-uniqueness}
There exists $\ka_1\in(0,\ka_2]$ such that, for every $0<\ka<\ka_1$,
the set $\cN_\ka$ has exactly one element.
\end{lemma}

\begin{proof}
Let $0<\ka<\ka_2$, $r:=r_\ka$, and $\rho,\mu\in\cN_\ka$, with rescaled
clusters $f_i(z)=r^d\rho_i(x_i+rz)$ and $g_i(z)=r^d\mu_i(x_i+rz)$, which are supported in $\overline{B_L}$. Set
\[
 \delta_i:=f_i-g_i,\qquad
 b_i:=\int\delta_i,\qquad
 p_i:=\int z\,\delta_i(z)\,dz.
\]
Write $b=(b_i)$ and $p=(p_i)$. On $\Omega_i:=\{f_i+g_i>0\}$ define
the mean diffusion coefficient
\begin{align*}
 q_i :=
  \frac{U'(f_i)-U'(g_i)}{f_i-g_i} = \frac{m}{m-1}\,\frac{f_i^{m-1}-g_i^{m-1}}{f_i-g_i},
\end{align*}
with $q_i:=mf_i^{m-2}$ where $f_i=g_i$. Note that $q_i > 0$ on $\Omega_i$.  Define
\[
 D:=\sum_i\int_{\Omega_i}q_i\delta_i^2\,dz,
 \qquad
 H:=\sum_i\|\delta_i\|_{L^1}.
\]
$D$ is finite since $\max(f_i,g_i) \le L$. Our goal is to show $D = 0$. 

To begin, $q_i$ is positively homogeneous of degree $m-2$ in $(f_i,g_i)$: with $s:=f_i+g_i$ and $\theta:=f_i/s$, we have $q_i=s^{m-2}\varphi(\theta)$ for a positive continuous function $\varphi$ on $[0,1]$. Hence $q_i\geq c_m(f_i+g_i)^{m-2}$ on $\Omega_i$ with $c_m=\min_{[0,1]}\varphi>0$. This implies $\int_{\Omega_i}q_i^{-1}\leq C$: If $1<m\leq2$, then $(f_i+g_i)^{2-m}\leq(2L)^{2-m}$. If $m>2$, then $(f_i+g_i)^{2-m}\leq \max(f_i,g_i)^{2-m} \le f_i^{2-m}\,\mathbf{1}_{\{f_i>0\}}+g_i^{2-m}\,\mathbf{1}_{\{g_i>0\}}$, and we use
\eqref{inverse-weight-bound}. Then the Cauchy--Schwarz inequality
$\|\delta_i\|_{L^1}^2\leq\int_{\Omega_i}q_i\delta_i^2\int_{\Omega_i}q_i^{-1}$ yields
\begin{equation}\label{difference-L1-control}
 H^2\leq C_1D.
\end{equation}

We now establish the inequality \eqref{direct-variational-inequality}. 
Both $\rho$ and $\mu$ are bounded, compactly supported, and absolutely
continuous, so $U'(\rho)$, $U'(\mu)$, $W*\rho$, and $W*\mu$ are bounded
on the common support and all integrals below are finite. \eqref{first_variation} gives
$\int(\ka U'(\rho)+W*\rho)\,d(\mu-\rho)\geq0$, and the same holds with
$\rho$ and $\mu$ exchanged.  Adding the two inequalities gives
\begin{equation}\label{ineq111}
 0\geq\ka\int\big(U'(\rho)-U'(\mu)\big)(\rho-\mu)\,dx
 +\iint W(x-y)\,d(\rho-\mu)(x)\,d(\rho-\mu)(y).
\end{equation}
Both measures are supported in $\bigcup_iB_{Lr}(x_i)$.  In each ball, the
change of variables $x=x_i+rz$, $\rho=r^{-d}f_i$, $\mu=r^{-d}g_i$,
$dx=r^d\,dz$, and $U'(s)=\frac{m}{m-1}s^{m-1}$ turn the first term into
\begin{equation}
 \ka r^{-d(m-1)}\sum_i\int_{\Omega_i}\frac{m}{m-1}
 \big(f_i^{m-1}-g_i^{m-1}\big)(f_i-g_i)\,dz
 =\ka r^{-d(m-1)}D=r^2D, \nonumber
\end{equation}
by \eqref{cluster-scale-balance}.  In the second term, for $x=x_i+rz$
and $y=x_j+rz'$ we have $x-y=e_{ij}+r(z-z')$, so by
\eqref{scaled-kernels}, $W(x-y)=w(1)+r^2K^r_{ij}(z,z')$ for $i\neq j$ and
$W(x-y)=r^2K^r_{ii}(z,z')$ for $i=j$.  The zeroth-order part of the
cross interactions is therefore
$w(1)\sum_{i\neq j}b_ib_j=w(1)\big((\sum_ib_i)^2-|b|^2\big)=c_W|b|^2$,
since $\sum_ib_i=0$. Hence \eqref{ineq111} becomes
\begin{equation}\label{direct-variational-inequality}
 0\geq c_W|b|^2+r^2(D+I),
 \qquad
 I:=\sum_i\iint K_{ii}^r\delta_i\delta_i
       +2\sum_{i<j}\iint K_{ij}^r\delta_i\delta_j,
\end{equation}
where each kernel is evaluated at $(z,z')$, and the two
density factors at $z$ and $z'$. 

Let $I^0$ denote $I$ with every kernel replaced by its limit in
\eqref{limiting-kernels}. We have
\[
 \iint|z-z'|^2\delta_i(z)\delta_i(z')\,dz\,dz'
 =2b_i\int|z|^2\delta_i\,dz-2|p_i|^2,
\]
\[
 \iint\big(e\cdot(z-z')\big)^2\delta_i(z)\delta_j(z')\,dz\,dz'
 =b_j\int(e\cdot z)^2\delta_i\,dz+b_i\int(e\cdot z)^2\delta_j\,dz
 -2(e\cdot p_i)(e\cdot p_j),
\]
for any unit vector $e$. Inserting these identities into $I^0$ yields $I^0=\cB(p)+R$, where
\[
 R= -\ibt\sum_ib_i\int|z|^2\delta_i\, dz
 +(\al-\bt)\sum_{i<j}\Big(
 b_j\int(e_{ij}\cdot z)^2\delta_i\, dz
 +b_i\int(e_{ij}\cdot z)^2\delta_j\, dz\Big).
\]
Since $\{ z : \delta_i(z) \ne 0\} \subset \overline{B_L}$,
$|R|\leq C|b|H$. Moreover we have $|I-I^0|\leq\ep(r)H^2$, where
\[
 \ep(r):=\max_{i,j}\sup_{\overline{B_L}\times\overline{B_L}}
 |K_{ij}^r-K_{ij}^0|\longrightarrow0
 \quad\hbox{as }r\downarrow0.
\]
Now in rescaled variables, the alignment conditions of $\rho$ read
\[
 \sum_i\int z\,f_i(z)\,dz=0,\qquad
 \sum_i\Big(\int z\,f_i(z)\,dz\Big)\cdot Sx_i=0\quad(S^T=-S),
\]
and the same identities hold for $g_i$.  Subtracting them shows that $p$
satisfies \eqref{moment-alignment}, and hence $\cB(p)\geq0$ by
Lemma \ref{dipole-positivity}. By \eqref{difference-L1-control}, we get the lower estimate
\begin{align*}
 D+I&\geq D-C|b|H-\ep(r)H^2\\
 &\geq\big(1-C_1\ep(r)\big)D-C|b|\sqrt{C_1D}
 \geq\Big(\frac12-C_1\ep(r)\Big)D-C_2|b|^2.
\end{align*}
Inserting this into \eqref{direct-variational-inequality} gives
\[
 0\geq\big(c_W-C_2r^2\big)|b|^2+r^2\Big(\frac12-C_1\ep(r)\Big)D.
\]
Choose $\ka_1\in(0,\ka_2]$ such that $C_2r_\ka^2\leq c_W/2$ and
$C_1\ep(r_\ka)\leq1/4$ for $\ka<\ka_1$.  Then both coefficients are
positive, so $D=0$.  Since $q_i>0$ on $\Omega_i$, this forces
$\delta_i=0$ for every $i$, and hence $\rho=\mu$. Finally, Lemma \ref{normalization-lemma} shows that $\cN_\ka$ is nonempty.
\end{proof}

\subsection{Proof of Theorem \ref*{main}}\label{subsection-main-proof}
We first normalize arbitrary minimizers, then
prove uniqueness, and finally use symmetry to identify the equal cluster
masses.

\begin{proof}[Proof of Theorem \ref{main}] Let $\ka_*:=\ka_1$ from Lemma~\ref{local-uniqueness}, and let $0<\ka<\ka_*$. By Lemma~\ref{existence}, $\cE_\ka$ has a global minimizer, and by Lemma~\ref{normalization-lemma} it is the image of an element of $\cN_\ka$ under a rigid motion. By Lemma~\ref{local-uniqueness}, $\cN_\ka$ consists of exactly one measure, which we denote by $\rho_\ka$. Hence every global minimizer is the image of $\rho_\ka$ under a rigid motion. Let $G$ be a symmetry of $\Delta^d$. By Lemma~\ref{normalization-lemma}, $G_\#\rho_\ka\in\cN_\ka$, hence $G_\#\rho_\ka=\rho_\ka$. For $i\neq j$, let $G_{ij}$ be the reflection in the hyperplane orthogonal to $x_i-x_j$. By \eqref{simplex-geometry}, $|x_i|=|x_j|$ and $x_k\cdot(x_i-x_j)=0$ for $k\notin\{i,j\}$, so $G_{ij}$ exchanges $x_i$ and $x_j$ and fixes the other vertices; thus it is a symmetry of $\Delta^d$. The barycenter $b$ of $\rho_\ka$ satisfies $G_{ij}b=b$, that is, $b\perp x_i-x_j$, for all $i\neq j$. These vectors span $\R^d$, so $b=0$. Since $\ka_*\leq\ka_3$ and $Lr_\ka<1/2$ for $\ka<\ka_3$, we have $Lr_{\ka_*}\leq1/2$; as $r_\ka$ is increasing in $\ka$, the balls $B_{Lr_\ka}(x_i)$ are pairwise disjoint for $0<\ka<\ka_*$. By Definition~\ref{normalized-definition}, $\supp\rho_\ka\subset\bigcup_iB_{Lr_\ka}(x_i)$. Since $G_{ij}$ maps $B_{Lr_\ka}(x_i)$ onto $B_{Lr_\ka}(x_j)$ and preserves $\rho_\ka$, the two balls carry the same mass. The masses of the $d+1$ disjoint balls sum to one, so each equals $1/(d+1)$. \end{proof}

\section{Proof of Theorem \ref*{shape}: Asymptotic shape of the clusters}\label{section-shape}
Throughout this section, $0<\ka<\ka_*$ and $\rho_\ka$ is the unique
element of $\cN_\ka$ from the proof of Theorem \ref{main}.  We write
$r=r_\ka$, $\rho_{\ka,i}$ for its clusters,
$f_{\ka,i}(z):=r^d\rho_{\ka,i}(x_i+rz)$ for its rescaled clusters,
$\ell_\ka$ for its multiplier, and $\Phi_\ka:=W*\rho_\ka$.  By
Theorem \ref{main}, $|\rho_{\ka,i}|=u$, and $G_\#\rho_\ka=\rho_\ka$ for
every symmetry $G$ of $\Delta^d$.

\begin{lemma}[The limiting cluster profile]\label{limiting-profile}
For each $i$, $f_{\ka,i}\to f_i$ in $L^1(\R^d)\cap L^\infty(\R^d)$ as
$\ka\to0$, where $f_i$ and $A_i$ are given by \eqref{intro-profile} and
\eqref{intro-matrix}, and $\lam_*>0$ is the unique constant satisfying
$\int f_i=u$ for every $i$. Moreover,
$\supp f_{\ka,i}\to E_i$ in Hausdorff distance.
\end{lemma}

\begin{proof}
Splitting $\Phi_\ka(x_i+rz)$ into
the contributions of the equal-mass clusters and using \eqref{scaled-kernels}
(the $d$ clusters $j\neq i$ each contribute $uw(1)$ at zeroth order) gives
\[
 \Phi_\ka(x_i+rz)=(1-u)w(1)+r^2\Psi_{\ka,i}(z),
 \qquad
 \Psi_{\ka,i}(z):=\sum_j\int K_{ij}^r(z,z')f_{\ka,j}(z')\,dz'.
\]
Recall $\ell_\Delta=(1-u)w(1)$, and put
$\lam_\ka:=(\ell_\ka-\ell_\Delta)/r^2$. By \eqref{multiplier-rate},
$|\lam_\ka|\leq C$.  By \eqref{eq123}, 
\begin{equation}\label{scaled-EL}
 f_{\ka,i}(z)=\left[\frac{m-1}{m}
 \big(\lam_\ka-\Psi_{\ka,i}(z)\big)\right]_+^{1/(m-1)},
 \qquad |z|<\frac{2\eta_0}{r}.
\end{equation}
$(f_{\ka,i})_i$ are bounded by $L$ and supported in
$\overline{B_L}$, and the convergence of the kernels gives local uniform $C^2$
bounds for $\Psi_{\ka,i}$, so
$f_{\ka,i}^{m-1}$ are uniformly Lipschitz on $\overline{B_{L+1}}$.  Let
$\ka_n\downarrow0$. After passing to a
subsequence, $\lam_{\ka_n}\to\bar\lam$, and $f_{\ka_n,i}^{m-1}$ converges uniformly on $\overline{B_{L+1}}$ by Arzel\`a--Ascoli, so $f_{\ka_n,i}$ also converges uniformly, and the common support bound gives convergence in $L^1(\R^d)$ and
$L^\infty(\R^d)$.  Write the limits as $f_i$. Then $\int f_i=u$, and
\[
 \Psi_{\ka_n,i}\longrightarrow
 \Psi_i(z):=\sum_j\int K_{ij}^0(z,z')f_j(z')\,dz'
 \quad\hbox{in }C^2\hbox{ on every bounded set}.
\]
Passing to the limit in \eqref{scaled-EL} yields
\begin{equation}\label{limiting-scaled-EL}
 f_i(z)=\left[\frac{m-1}{m}
 \big(\bar\lam-\Psi_i(z)\big)\right]_+^{1/(m-1)}
 \quad\hbox{for all }z\in\R^d.
\end{equation}
At this stage, the limits may depend on the chosen subsequence. We will show that $\bar\lam-\Psi_i(z)=\lam_*-z\cdot A_i z$, independently of this choice.

Let $a_i:=u^{-1}\int zf_i(z)\,dz$. We will show $a_i = 0$. If $G$ is a symmetry with
$Gx_i=x_l$, then $G$ maps $B_{1/2}(x_i)$ onto $B_{1/2}(x_l)$, and
$G_\#\rho_\ka=\rho_\ka$ gives $f_{\ka,l}(Gz)=f_{\ka,i}(z)$; hence
$a_l=Ga_i$.  For $j,k\neq i$, the reflection in the hyperplane
orthogonal to $x_j-x_k$ fixes $x_i$, so it fixes $a_i$, hence $a_i \cdot(x_j -x_k) = 0$. Since the vectors
$x_j-x_k$ with $j,k\neq i$ span $x_i^\perp$, $a_i=t_ix_i$ for some
$t_i\in\R$, and $a_l=Ga_i$ shows that
$t_i=t$ does not depend on $i$.

Translate just the $i$-th cluster of $\rho_\ka$ by a small vector $v$,
leaving the other clusters fixed, and let $\rho^v$ denote the resulting measure.  Their supports remain disjoint. Since the internal energy and the
self-interaction of the translated cluster are unchanged,
\[
 \cE_\ka(\rho^v)-\cE_\ka(\rho_\ka)
 =\sum_{j\neq i}\iint\big(W(x+v-y)-W(x-y)\big)
   \,d\rho_{\ka,i}(x)\,d\rho_{\ka,j}(y).
\]
Minimality at $v=0$ in all
directions gives, by differentiation in $v$, the vector identity
\begin{equation}\label{cluster-force-balance}
 \sum_{j\neq i}\iint\grad W(x-y)
 \,d\rho_{\ka,i}(x)\,d\rho_{\ka,j}(y)=0.
\end{equation}
Since $\grad W(e_{ij})=0$ and
$D^2W(e_{ij})=(\al-\bt)e_{ij}\otimes e_{ij}$, Taylor's formula gives
\[
 \frac{\grad W(e_{ij}+r(z-z'))}{r}
 \longrightarrow
 (\al-\bt)(e_{ij}\otimes e_{ij})(z-z')
\]
uniformly for $z,z'\in\overline{B_L}$. We rescale
\eqref{cluster-force-balance}, divide by $r$, and pass to the limit
along the subsequence. Using the $L^1$ convergence of the rescaled clusters, the limit is
\[
 0 = (\al-\bt)\sum_{j\neq i}(e_{ij}\otimes e_{ij})
 \iint(z-z')f_i(z)f_j(z')\,dz\,dz'
 =(\al-\bt)u^2\sum_{j\neq i}(e_{ij}\otimes e_{ij})(a_i-a_j).
\]
Hence $\sum_{j\neq i}(e_{ij}\otimes e_{ij})(a_i-a_j)=0$. 
With $a_i=tx_i$, $(e_{ij}\otimes e_{ij})(a_i-a_j)=t\,e_{ij}|e_{ij}|^2=t\,e_{ij}$,
hence $t\sum_{j\neq i}e_{ij}=t(d+1)x_i = 0$. Hence
$t=0$, and all limiting cluster means $a_i$ vanish.

Since each limiting cluster has mass $u$ and zero first moment, direct computation gives $\nabla \Psi_i (0) = 0$ and $\nabla^2\Psi_i(z) = u \big( (\al- \bt)  \sum_{j \ne i} (e_{ij} \otimes e_{ij}) - \ibt I_d\big)$. Lemma \ref{simplex-wells} with \eqref{simplex-well-hessian} gives $\nabla^2\Psi_i(z) =  \grad^2V_\Del(x_i) = 2A_i$. As $\Psi_i$ is quadratic, $\Psi_i(z)= \Psi_i(0) + z\cdot A_iz$. Hence \eqref{limiting-scaled-EL} is \eqref{intro-profile} with $\lam_*$ replaced by $\lam_i:=\bar\lam-\Psi_i(0)$. Since $A_i$ is positive
definite, the positive mass forces $\lam_i>0$.  The change of
variables $z=\sqrt{\lam_i}\,A_i^{-1/2}y$ gives
\begin{equation}\label{profile-mass}
 \int f_i=\Big(\frac{m-1}{m}\Big)^{\frac1{m-1}}\omega_{d,m}
 (\det A_i)^{-1/2}\lam_i^{\frac{d(m-1)+2}{2(m-1)}},
 \quad
 \omega_{d,m}:=\int_{B_1}(1-|y|^2)^{\frac1{m-1}}\,dy,
\end{equation}
which increases strictly from zero to infinity with $\lam_i$.  By
\eqref{intro-matrix}, $\det A_i=(u/2)^da_\parallel a_\perp^{d-1}$ does
not depend on $i$ (for $d=1$ the factor $a_\perp$ is absent).  Hence the
mass condition $\int f_i = u$ determines a common value $\lam_i=\lam_*$, and this determines $f_i$ by \eqref{intro-profile}.  Since the limit does not depend on the subsequence, the whole family converges as $\ka\to0$. 

By the same argument, $\lam_\ka-\Psi_{\ka,i}\to\lam_*-z\cdot A_iz$ uniformly on $\overline{B_{L+1}}$, and $\supp f_{\ka,i}$ is the closure of the set $\{ z \in \overline{B_L} : \lam_\ka-\Psi_{\ka,i}(z) > 0 \}$. 
Let $E_i^{+} = \{z\in\overline{B_{L+1}}:\operatorname{dist}(z,E_i)\leq\ep\}$ and $E_i^{-} = \{z\in E_i :\operatorname{dist}(z, E_i^c)\geq\ep\}$ given $\ep >0$. The limit $\lam_*-z\cdot A_iz$ is negative and bounded away from zero on $\overline{B_{L+1}} \setminus E_i^{+}$, while positive and bounded away from zero on $E_i^{-}$. Hence the uniform convergence  gives Hausdorff convergence of $\supp f_{\ka,i}$ to $E_i$.
\end{proof}

\begin{remark}[Explicit constants]\label{explicit-constants}
Since $\omega_{d,m}=\pi^{d/2}\Gamma\big(\frac m{m-1}\big)/
\Gamma\big(\frac m{m-1}+\frac d2\big)$, \eqref{profile-mass} gives
\begin{equation}\label{lambda-explicit}
 \lam_*=\left[\frac{u\,(u/2)^{d/2}\,a_\parallel^{1/2}\,a_\perp^{(d-1)/2}}
 {\big(\frac{m-1}{m}\big)^{1/(m-1)}\,\omega_{d,m}}
 \right]^{\frac{2(m-1)}{d(m-1)+2}},
\end{equation}
where the factor $a_\perp$ is absent for $d=1$. For the Wasserstein asymptotics below, we also record that \eqref{intro-matrix} gives \[ \tr(A_i^{-1}) =\frac{2}{u}\left(\frac{1}{a_\parallel} +\frac{d-1}{a_\perp}\right) \qquad(d\geq2), \] while $\tr(A_i^{-1})=2/(u a_\parallel)$ when $d=1$.
\end{remark}

\begin{proof}[Proof of Theorem \ref{shape}]
Fix $0<\ka<\ka_*$ and write $r:=r_\ka$.  Since
$\rho_\ka\in\cN_\ka$, \eqref{normalized-estimates} gives
$\|\rho_{\ka,i}\|_{L^\infty}\leq Lr^{-d}$ for every $i$. Let $O_{\ka,i}:=\{\rho_{\ka,i}>0\}$.  Since $\ka_*\leq\ka_2$, Lemma~\ref{inverse-weight-lemma} gives
$D^2\Phi_\ka\geq c_*I_d$ on $B_{2\eta_0}(x_i)$, and the choice of
$\ka_2$ ensures $Lr\leq\eta_0$.  Thus $\supp\rho_{\ka,i}\subset B_{Lr}(x_i)\subset B_{\eta_0}(x_i)$. By \eqref{rho-version} and the separation of the clusters,
\[
 O_{\ka,i}
 =\{x\in B_{2\eta_0}(x_i):\Phi_\ka(x)<\ell_\ka\}.
\]
Hence $O_{\ka,i}$ is a nonempty convex open set 
contained in $B_{\eta_0}(x_i)$, and since $\Phi_\ka$ is convex, $\partial O_{\ka,i} =\{x\in B_{2\eta_0}(x_i):\Phi_\ka(x)=\ell_\ka\}$. 
Moreover $\grad\Phi_\ka\neq0$ on $\partial O_{\ka,i}$: by strict convexity, a critical point would minimize $\Phi_\ka$ on $B_{2\eta_0}(x_i)$, contradicting $\Phi_\ka<\ell_\ka$ on $O_{\ka,i}$. 
Since $\Phi_\ka\in C^2$, the implicit function theorem
shows that $\partial O_{\ka,i}$ is $C^2$, and on $O_{\ka,i}$, $\rho_\ka^{m-1}
 =\frac{m-1}{m\ka}(\ell_\ka-\Phi_\ka)$, which is concave. Now Lemma~\ref{limiting-profile} gives the convergence of the rescaled
densities in $L^1(\R^d)\cap L^\infty(\R^d)$ and the Hausdorff
convergence of their supports.  The common normalization constant
$\lam_*$ is given by \eqref{lambda-explicit}.

It remains to prove the Wasserstein asymptotics.  Since $Lr<1/2$
and the vertices have mutual distance one, every point in the
$i$-th cluster is closer to $x_i$ than to any other vertex.
By Theorem~\ref{main}, each cluster has mass $u$, so the transport plan sending the entire $i$-th cluster to $x_i$ is optimal for $d_2$ and has second marginal $\rho_\Delta$. Hence
\[
 d_2(\rho_\ka,\rho_\Delta)^2
 =\sum_{i=0}^d\int|x-x_i|^2\,d\rho_{\ka,i}(x)
 =r^2\sum_{i=0}^d\int|z|^2f_{\ka,i}(z)\,dz.
\]
The $L^1$ convergence in Lemma~\ref{limiting-profile} and the common
support bound $\supp f_{\ka,i}\subset\overline{B_L}$ give
\[
 \lim_{\ka\to0}\frac{d_2(\rho_\ka,\rho_\Delta)^2}{r_\ka^2}
 =\sum_{i=0}^d\int|z|^2f_i(z)\,dz.
\]
Put $s=\frac{1}{m-1}$.  The change of variables
$z=\sqrt{\lam_*}\,A_i^{-1/2}y$, together with $u = \int f_i$, yields
\[
 \frac{\int|z|^2f_i(z)\,dz}{u}
 =\lam_*
 \frac{\displaystyle\int_{B_1}
       y\cdot A_i^{-1}y\,(1-|y|^2)^s\,dy}
      {\displaystyle\int_{B_1}(1-|y|^2)^s\,dy}.
\]
By radial symmetry of the weight,
\[
 \int_{B_1}y\cdot A_i^{-1}y\,(1-|y|^2)^s\,dy
 =\frac{\tr(A_i^{-1})}{d}
   \int_{B_1}|y|^2(1-|y|^2)^s\,dy.
\]
The Beta-function identity
$\displaystyle 
 \frac{\int_{B_1}|y|^2(1-|y|^2)^s\,dy}
      {\int_{B_1}(1-|y|^2)^s\,dy}
 =\frac{d}{d+2s+2}
$
thus gives
\[
 \int|z|^2f_i(z)\,dz
 =\frac{u\lam_*\tr(A_i^{-1})}{d+2+2/(m-1)}.
\]
Summing over $i$, using $(d+1)u=1$, and substituting 
$\tr(A_i^{-1})
 =\frac{2}{u}\left(\frac{1}{a_\parallel}
                  +\frac{d-1}{a_\perp}\right)$, 
with the second term omitted when $d=1$, proves the asserted limit.
\end{proof}

\section{Concluding remarks}\label{section-remarks}

Small porous-medium diffusion replaces the atomic simplex ground state
by diffuse clusters while preserving its uniqueness up to rigid
motions; because uniqueness is proved before symmetry is invoked, the
exact equality of the cluster masses is a consequence rather than an
assumption.  The boundary of the strict simplex regime behaves
differently.  For $d\geq2$ we have $\al_\Delta(2)=4$ \cite{DLM23}, so
$(\al,\bt)=(4,2)\notin\cS$.  At this pair $\rho_\Delta$ is still a
minimizer of $\cW$ \cite{DLM23}, but small diffusion selects a different
geometry.

\begin{proposition}[Selection at the degenerate pair $(4,2)$]
\label{critical-selection}
Let $d\geq2$, $m>1$, $(\al,\bt)=(4,2)$, and 
$R=\sqrt{d/(2d+2)}$.  For every $\ka>0$, $\cE_\ka$ has exactly one global minimizer $\rho_\ka$ in $\cP_0(\R^d)$; it is radially symmetric, and every global minimizer of $\cE_\ka$ is a translate of $\rho_\ka$.  As $\ka\downarrow0$, $\rho_\ka$ converges
narrowly to the normalized surface measure $\sigma_R$ on
$\partial B_R$.
\end{proposition}

\begin{proof}
Write $\cW_{4,2}$ for the interaction energy at this pair.  Measures of
finite energy have finite fourth moments by
\eqref{centered-moment-bound}.  For such $\rho\in\cP_0(\R^d)$, put
$C_\rho:=\int x\otimes x\,d\rho(x)$ and $s_\rho:=\tr C_\rho$.
Expanding $|x-y|^4=(|x|^2+|y|^2-2x\cdot y)^2$ and using
$\int x\,d\rho(x)=0$, we get
\begin{equation*}
 \cW_{4,2}(\rho)
 =\frac14\int|x|^4\,d\rho(x)
  +\frac14s_\rho^2+\frac12|C_\rho|^2-\frac12s_\rho.
\end{equation*}
Along linear interpolations in $\cP_0(\R^d)$, the quantities $\int|x|^4\,d\rho$, $s_\rho$, and $C_\rho$ are affine, so the right-hand side is convex; this convexity argument goes back to Lopes \cite{Lop19} and was used in \cite[Lemma~2.2]{DLM23}. Since $U$ is strictly convex, $\cE_\ka$ is strictly convex along linear interpolations of centered densities of finite energy.

The proof of Lemma~\ref{existence} applies at this pair, so $\cE_\ka$
has a global minimizer, which we may translate into $\cP_0(\R^d)$.
Strict convexity gives uniqueness in $\cP_0(\R^d)$; hence every global
minimizer, translated by its barycenter, equals $\rho_\ka$.  For
$Q\in O(d)$, $Q_\#\rho_\ka$ is again a global minimizer in
$\cP_0(\R^d)$, so $Q_\#\rho_\ka=\rho_\ka$. The same proof gives a fixed ball containing $\supp\rho_\ka$ for
$0<\ka\leq1$.  Hence every sequence $\ka_n\downarrow0$ has a
subsequence along which $\rho_{\ka_n}$ converges narrowly to an
orthogonally invariant $\rho_0\in\cP_0(\R^d)$, with
$\cW_{4,2}(\rho_{\ka_n})\to\cW_{4,2}(\rho_0)$ because $W$ is continuous
and the supports are uniformly bounded. For every smooth compactly
supported probability density $\nu$,
\[
 \cW_{4,2}(\rho_{\ka_n})\leq\cE_{\ka_n}(\rho_{\ka_n})
 \leq\cE_{\ka_n}(\nu)\longrightarrow\cW_{4,2}(\nu).
\]
Taking for $\nu$ mollifications of $\sigma_R$ gives
$\cW_{4,2}(\rho_0)\leq\cW_{4,2}(\sigma_R)=\min\cW_{4,2}$, where the
equality holds by \cite[Theorem 1.1]{DLM23}.  Thus $\rho_0$ is a
centered minimizer of $\cW_{4,2}$, and by the same theorem it is
supported on $\partial B_R$.  Being orthogonally invariant, it equals
$\sigma_R$.  Since the limit does not depend on the subsequence, the
whole family converges.
\end{proof}

Proposition~\ref{critical-selection} shows that our estimates cannot
be uniform as $\al\downarrow4$ when $\bt=2$: the tangential eigenvalue
$ua_\perp=u(\al-4)/2$ of the simplex wells degenerates, and the
selected geometry changes from separated clusters to a uniform shell.
The estimates are likewise not uniform as $\al\to\infty$ or
$m\to\infty$; these limits lead to the isodiametric problem of
\cite{LM21} and the density-constrained problem of \cite{CT20}.  

For other radial kernels $W(x)=w(|x|)$ whose zero-diffusion minimizer is
the equally weighted simplex, the localization and normalization steps
carry over once the vertices are isolated, non-degenerate minima of the
simplex potential.  The uniqueness step uses in addition that
$w'(1)=0$, $w''(1)>0$, and $D^2W(0)\leq0$; the limiting interaction then
has the form of Lemma~\ref{dipole-positivity}, with the same positivity
on aligned first moments.

Several questions remain open: whether uniqueness up to rigid motions
persists for all $\ka>0$; how the clusters merge into the shell in the
joint limit $\ka\downarrow0$, $\al\downarrow4$ with $\bt=2$; whether, for
linear diffusion, the clusters rescaled at the scale $\ka^{1/2}$
converge to Gaussian profiles proportional to $e^{-z\cdot A_iz}$, the
$m\downarrow1$ limits of \eqref{intro-profile}; and whether, for
$\al>\bt>2$, the unequal-mass simplex configurations that are local
minimizers of $\cW$ \cite{LM21} give rise to steady states or local
minimizers of $\cE_\ka$ with unequal cluster masses.


\bibliographystyle{plain} 
\bibliography{LimSeo_ground_states}

@article{FM25,
	author = {Frank, Rupert L. and Matzke, Ryan W.},
	date = {2025/02/11},
	doi = {10.1007/s00205-025-02084-1},
	id = {Frank2025},
	isbn = {1432-0673},
	journal = {Archive for Rational Mechanics and Analysis},
	number = {2},
	pages = {15},
	title = {Minimizers for an Aggregation Model with Attractive--Repulsive Interaction},
	url = {https://doi.org/10.1007/s00205-025-02084-1},
	volume = {249},
	year = {2025}}

@book{AGS08,
	address = {Basel},
	author = {Ambrosio, Luigi and Gigli, Nicola and Savar{\'e}, Giuseppe},
	edition = {Second},
	publisher = {Birkh{\"a}user Verlag},
	series = {Lectures in Mathematics ETH Z{\"u}rich},
	title = {Gradient Flows in Metric Spaces and in the Space of Probability Measures},
	year = {2008}}

@article{BCLR13,
	author = {Balagu{\'e}, Daniel and Carrillo, Jos{\'e} A. and Laurent, Thomas and Raoul, Ga{\"e}l},
	doi = {10.1007/s00205-013-0644-6},
	journal = {Arch. Ration. Mech. Anal.},
	number = {3},
	pages = {1055--1088},
	title = {Dimensionality of local minimizers of the interaction energy},
	volume = {209},
	year = {2013}}

@article{BD08,
	author = {Burger, Martin and Di Francesco, Marco},
	doi = {10.3934/nhm.2008.3.749},
	journal = {Netw. Heterog. Media},
	number = {4},
	pages = {749--785},
	title = {Large time behavior of nonlocal aggregation models with nonlinear diffusion},
	volume = {3},
	year = {2008}}

@article{BDF13,
	author = {Burger, Martin and Di Francesco, Marco and Franek, Marzena},
	doi = {10.4310/CMS.2013.v11.n3.a3},
	journal = {Commun. Math. Sci.},
	number = {3},
	pages = {709--738},
	title = {Stationary states of quadratic diffusion equations with long-range attraction},
	volume = {11},
	year = {2013}}

@article{BFH14,
	author = {Burger, Martin and Fetecau, Razvan C. and Huang, Yanghong},
	doi = {10.1137/130923786},
	journal = {SIAM J. Appl. Dyn. Syst.},
	number = {1},
	pages = {397--424},
	title = {Stationary states and asymptotic behavior of aggregation models with nonlinear local repulsion},
	volume = {13},
	year = {2014}}

@article{CCP15,
	author = {Ca{\~n}izo, Jos{\'e} A. and Carrillo, Jos{\'e} A. and Patacchini, Francesco S.},
	doi = {10.1007/s00205-015-0852-3},
	journal = {Arch. Ration. Mech. Anal.},
	number = {3},
	pages = {1197--1217},
	title = {Existence of compactly supported global minimisers for the interaction energy},
	volume = {217},
	year = {2015}}

@incollection{CCY19,
	address = {Cham},
	author = {Carrillo, Jos{\'e} A. and Craig, Katy and Yao, Yao},
	booktitle = {Active Particles, Volume 2},
	doi = {10.1007/978-3-030-20297-2_3},
	pages = {65--108},
	publisher = {Birkh{\"a}user/Springer},
	series = {Model. Simul. Sci. Eng. Technol.},
	title = {Aggregation-diffusion equations: dynamics, asymptotics, and singular limits},
	year = {2019}}

@article{CDM16,
	author = {Carrillo, Jos{\'e} A. and Delgadino, Matias G. and Mellet, Antoine},
	doi = {10.1007/s00220-016-2598-7},
	journal = {Comm. Math. Phys.},
	number = {3},
	pages = {747--781},
	title = {Regularity of local minimizers of the interaction energy via obstacle problems},
	volume = {343},
	year = {2016}}

@article{CDP19,
	author = {Carrillo, Jos{\'e} A. and Delgadino, Matias G. and Patacchini, Francesco S.},
	doi = {10.1142/S0219530518500276},
	journal = {Anal. Appl. (Singap.)},
	number = {3},
	pages = {393--423},
	title = {Existence of ground states for aggregation-diffusion equations},
	volume = {17},
	year = {2019}}

@article{CFP17,
	author = {Carrillo, Jos{\'e} A. and Figalli, Alessio and Patacchini, Francesco S.},
	doi = {10.1016/j.anihpc.2016.10.004},
	journal = {Ann. Inst. H. Poincar{\'e} C Anal. Non Lin{\'e}aire},
	number = {5},
	pages = {1299--1308},
	title = {Geometry of minimizers for the interaction energy with mildly repulsive potentials},
	volume = {34},
	year = {2017}}

@article{CFT15,
	author = {Choksi, Rustum and Fetecau, Razvan C. and Topaloglu, Ihsan},
	doi = {10.1016/j.anihpc.2014.09.004},
	journal = {Ann. Inst. H. Poincar{\'e} C Anal. Non Lin{\'e}aire},
	number = {6},
	pages = {1283--1305},
	title = {On minimizers of interaction functionals with competing attractive and repulsive potentials},
	volume = {32},
	year = {2015}}

@article{CHVY19,
	author = {Carrillo, Jos{\'e} A. and Hittmeir, Sabine and Volzone, Bruno and Yao, Yao},
	doi = {10.1007/s00222-019-00898-x},
	journal = {Invent. Math.},
	number = {3},
	pages = {889--977},
	title = {Nonlinear aggregation-diffusion equations: radial symmetry and long time asymptotics},
	volume = {218},
	year = {2019}}

@article{CMV03,
	author = {Carrillo, Jos{\'e} A. and McCann, Robert J. and Villani, C{\'e}dric},
	doi = {10.4171/RMI/376},
	journal = {Rev. Mat. Iberoamericana},
	number = {3},
	pages = {971--1018},
	title = {Kinetic equilibration rates for granular media and related equations: entropy dissipation and mass transportation estimates},
	volume = {19},
	year = {2003}}

@article{CMV06,
	author = {Carrillo, Jos{\'e} A. and McCann, Robert J. and Villani, C{\'e}dric},
	doi = {10.1007/s00205-005-0386-1},
	journal = {Arch. Ration. Mech. Anal.},
	number = {2},
	pages = {217--263},
	title = {Contractions in the 2-{W}asserstein length space and thermalization of granular media},
	volume = {179},
	year = {2006}}

@article{CT20,
	author = {Craig, Katy and Topaloglu, Ihsan},
	doi = {10.1016/j.anihpc.2019.10.003},
	journal = {Ann. Inst. H. Poincar{\'e} C Anal. Non Lin{\'e}aire},
	number = {2},
	pages = {239--279},
	title = {Aggregation-diffusion to constrained interaction: minimizers \& gradient flows in the slow diffusion limit},
	volume = {37},
	year = {2020}}

@article{DFJ19,
	author = {Di Francesco, Marco and Jaafra, Yahya},
	doi = {10.3934/krm.2019013},
	journal = {Kinet. Relat. Models},
	number = {2},
	pages = {303--322},
	title = {Multiple large-time behavior of nonlocal interaction equations with quadratic diffusion},
	volume = {12},
	year = {2019}}

@article{DLMshell22,
	author = {Davies, Cameron and Lim, Tongseok and McCann, Robert J.},
	doi = {10.1137/21M1455309},
	journal = {SIAM J. Appl. Math.},
	number = {4},
	pages = {1520--1536},
	title = {Classifying minimum energy states for interacting particles: spherical shells},
	volume = {82},
	year = {2022}}

@article{DLM23,
	author = {Davies, Cameron and Lim, Tongseok and McCann, Robert J.},
	doi = {10.1007/s00220-022-04564-x},
	journal = {Comm. Math. Phys.},
	number = {2},
	pages = {577--598},
	title = {Classifying minimum energy states for interacting particles: regular simplices},
	volume = {399},
	year = {2023}}

@article{DYY22,
	author = {Delgadino, Matias G. and Yan, Xukai and Yao, Yao},
	doi = {10.1002/cpa.21950},
	journal = {Comm. Pure Appl. Math.},
	number = {1},
	pages = {3--59},
	title = {Uniqueness and nonuniqueness of steady states of aggregation-diffusion equations},
	volume = {75},
	year = {2022}}

@article{FJ21,
	author = {Fagioli, Simone and Jaafra, Yahya},
	doi = {10.3934/nhm.2021010},
	journal = {Netw. Heterog. Media},
	number = {3},
	pages = {377--411},
	title = {Multiple patterns formation for an aggregation/diffusion predator-prey system},
	volume = {16},
	year = {2021}}

@article{FKT19,
	author = {Fetecau, Razvan C. and Kovacic, Mitchell and Topaloglu, Ihsan},
	doi = {10.3934/dcdsb.2018238},
	journal = {Discrete Contin. Dyn. Syst. Ser. B},
	number = {4},
	pages = {1815--1842},
	title = {Swarming in domains with boundaries: approximation and regularization by nonlinear diffusion},
	volume = {24},
	year = {2019}}

@article{Kai17,
	author = {Kaib, Gunnar},
	doi = {10.1137/16M1072450},
	journal = {SIAM J. Math. Anal.},
	number = {1},
	pages = {272--296},
	title = {Stationary states of an aggregation equation with degenerate diffusion and bounded attractive potential},
	volume = {49},
	year = {2017}}

@article{KKLS21,
	author = {Kang, Kyungkeun and Kim, Hwa Kil and Lim, Tongseok and Seo, Geuntaek},
	doi = {10.1007/s00526-020-01882-7},
	journal = {Calc. Var. Partial Differential Equations},
	number = {1},
	pages = {Paper No. 15},
	title = {Uniqueness and characterization of local minimizers for the interaction energy with mildly repulsive potentials},
	volume = {60},
	year = {2021}}

@article{LM21,
	author = {Lim, Tongseok and McCann, Robert J.},
	doi = {10.1007/s00205-021-01632-9},
	journal = {Arch. Ration. Mech. Anal.},
	number = {2},
	pages = {553--576},
	title = {Isodiametry, variance, and regular simplices from particle interactions},
	volume = {241},
	year = {2021}}

@unpublished{LS-dyn,
	author = {Lim, Tongseok and Seo, Geuntaek},
	note = {Companion paper, in preparation},
	title = {Local asymptotic stability of simplex ground states for aggregation--diffusion equations with mild repulsion},
	year = {2026}}

@article{Lop19,
	author = {Lopes, Orlando},
	doi = {10.3934/cpaa.2019102},
	journal = {Commun. Pure Appl. Anal.},
	number = {5},
	pages = {2265--2282},
	title = {Uniqueness and radial symmetry of minimizers for a nonlocal variational problem},
	volume = {18},
	year = {2019}}

@book{San15,
	address = {Cham},
	author = {Santambrogio, Filippo},
	publisher = {Birkh{\"a}user/Springer},
	series = {Progress in Nonlinear Differential Equations and their Applications},
	title = {Optimal Transport for Applied Mathematicians},
	volume = {87},
	year = {2015}}

@article{SST15,
	author = {Simione, Robert and Slep{\v{c}}ev, Dejan and Topaloglu, Ihsan},
	doi = {10.1007/s10955-015-1215-z},
	journal = {J. Stat. Phys.},
	number = {4},
	pages = {972--986},
	title = {Existence of ground states of nonlocal-interaction energies},
	volume = {159},
	year = {2015}}

@book{V07,
	address = {Oxford},
	author = {V{\'a}zquez, Juan Luis},
	doi = {10.1093/acprof:oso/9780198569039.001.0001},
	publisher = {The Clarendon Press, Oxford University Press},
	series = {Oxford Mathematical Monographs},
	title = {The Porous Medium Equation: Mathematical Theory},
	year = {2007}}

\end{document}